\documentclass[11pt,reqno]{amsart}
\pdfoutput=1
\usepackage{amsthm, amsmath,amsfonts,amssymb,euscript,hyperref,graphics,color,slashed}
\usepackage{graphicx}
\usepackage{mathrsfs}
\usepackage{comment}
\usepackage{import}
\usepackage{tikz}
\usepackage{latexsym}
\usepackage[makeroom]{cancel}
\usepackage[font=normalsize]{caption}
\usepackage{subfigure}
\usepackage[title]{appendix}
\usepackage{enumerate}
\usepackage{mathrsfs}
\usepackage{empheq}   
\usepackage{bm}
\usepackage{bookmark}
\usepackage{arydshln}
\usepackage{lineno,hyperref}
\usepackage{cite}
\usepackage{float}
\usepackage[utf8]{inputenc}
\usepackage{dashbox}
\usepackage{ulem}
\usepackage{lipsum}

\def\inte#1{
\displaystyle\mathop{#1\kern0pt}^\circ }

\def\virgp{\raise 2pt\hbox{,}}
\def\cdotpv{\raise 2pt\hbox{;}}

\def\C{\mathop{\mathbb C\kern 0pt}\nolimits}
\def\DD{\mathop{\mathbb D\kern 0pt}\nolimits}
\def\EE{\mathop{{\mathbb E \kern 0pt}}\nolimits}
\def\K{\mathop{\mathbb K\kern 0pt}\nolimits}
\def\N{\mathop{\mathbb N\kern 0pt}\nolimits}
\def\Q{\mathop{\mathbb Q\kern 0pt}\nolimits}
\def\R{\mathop{\mathbb R\kern 0pt}\nolimits}
\def\SS{\mathop{\mathbb S\kern 0pt}\nolimits}
\def\ZZ{\mathop{\mathbb Z\kern 0pt}\nolimits}
\def\TT{\mathop{\mathbb T\kern 0pt}\nolimits}
\def\P{\mathop{\mathbb P\kern 0pt}\nolimits}

\newcommand{\beq}{\begin{equation}}
\newcommand{\eeq}{\end{equation}}
\newcommand{\ben}{\begin{eqnarray}}
\newcommand{\een}{\end{eqnarray}}
\newcommand{\beno}{\begin{eqnarray*}}
\newcommand{\eeno}{\end{eqnarray*}}

\newtheorem{thm}{Theorem}[section]
\newtheorem{lem}{Lemma}[section]
\newtheorem{rmk}{Remark}[section]

\newtheorem{prop}{Proposition}[section]

\newtheorem*{Main Theorem}{Main Theorem}

\numberwithin{equation}{section}

\allowdisplaybreaks

\title[$H^{\frac{11}{4}}(\mathbb{R}^2)$ ill-posedness of 2D Elastic Waves]{$H^{\frac{11}{4}}(\mathbb{R}^2)$ ill-posedness for 2D Elastic Wave system}

\author{Xinliang An$^*$$^1$}\author{Haoyang Chen$^*$$^2$}\author{Silu Yin$^{\dag}$$^3$}

\date{}

\AtEndDocument{%
  \par
  \medskip
  \begin{tabular}{@{}l@{}}%
    \textsc{$^*$Department of Mathematics, National University of Singapore, Singapore}\\
    \textit{E-mail address}: \texttt{$^1$matax@nus.edu.sg}\\
    \textit{E-mail address}: \texttt{$^2$hychen@nus.edu.sg}\\
    \textsc{$^\dag$School of Mathematics, Hangzhou Normal University, China}\\
    \textit{E-mail address}: \texttt{$^3$yins@hznu.edu.cn}\\
  \end{tabular}}

\begin{document}

\maketitle
\newcommand\blfootnote[1]{%
\begingroup
\renewcommand\thefootnote{}\footnote{#1}%
\addtocounter{footnote}{-1}%
\endgroup
}
\begin{abstract}
In this paper, we prove that for the 2D elastic wave equations, a physical system with multiple wave-speeds, its Cauchy problem fails to be locally well-posed in $H^{\frac{11}{4}}(\mathbb{R}^2)$. The ill-posedness here is driven by instantaneous shock formation. In 2D Smith-Tataru showed that the Cauchy problem for a single quasilinear wave equation is locally well-posed in $H^s$ with $s>\frac{11}{4}$. Hence our $H^{\frac{11}{4}}$ ill-posedness obtained here is a desired result. Our proof relies on combining a geometric method and an algebraic wave-decomposition approach, equipped with detailed analysis of the corresponding hyperbolic system.
\end{abstract}

\blfootnote{XA is supported by NUS startup grant R-146-000-269-133, and MOE Tier 1 grants R-146-000-321-144 and R-146-000-335-114. HC acknowledges the support of MOE Tier 1 grant R-146-000-335-114 and NSFC (Grant No. 12171097). SY is supported by NSFC (Grant No. 12001149) and Zhejiang Provincial Natural Science Foundation of China (Grant No. LQ19A010006).}
\begin{center}
\textit{To Professor Demetrios Christodoulou with admiration and gratitude.
}
\end{center}
\tableofcontents

\section{Introduction}
This article investigates the low regularity ill-posedness for elastic waves in two spatial dimensions. For homogeneous isotropic hyperelastic materials, the motion of displacement $U=(U_1,U_2)$ satisfies a quasilinear wave system with multiple wave-speeds:
\begin{equation} \label{elaswave}
    \partial_t^2 U-c_2^2\Delta U-(c_1^2-c_2^2)\nabla(\nabla\cdot U)=G(\nabla U,\nabla^2U),
\end{equation}
where $c_1,c_2$ are two constants satisfying $c_1>c_2>0$. The precise form of $G(\nabla U,\nabla^2U)$ is discussed in Section \ref{eq and main thm}.

The mathematical study of the elastic wave system was pioneered by Fritz John. For the 3D elastic wave equations, John proved in \cite{john84} that the singularities could arise from (smooth) small initial data with radial symmetry. When the elastic wave system satisfies null conditions, Agemi \cite{Agemi} and Sideris \cite{Sideris96} proved that the global existence of solutions to the Cauchy problem of \eqref{elaswave} holds for small initial data\footnote{For 2D elastic wave equations, under null conditions, Zha \cite{zha} showed the global existence of classical solutions for sufficiently small radially-symmetric initial data.}. In \cite{an}, the $H^3$ ill-posedness of the 3D elastic wave system was proved by An-Chen-Yin. In this paper, we focus on the 2D elastic wave equations. The main conclusion we obtain is that \textit{the Cauchy problem for the 2D elastic wave equations is ill-posed in $H^{\frac{11}{4}}(\mathbb{R}^2)$.} This $H^\frac{11}{4}$ ill-posedness is consistent with the low-regularity local well-posedness result for the quasilinear scalar wave equation by Smith-Tataru \cite{tataru}. In this paper, we will further show that the mechanism driving the ill-posedness is the shock formation, which is characterized by the inverse density for the characteristics being zero.

Our research is motivated by a series of classic works on low-regularity ill-posedness and shock formation. Under planar symmetry, Lindblad \cite{lindblad93,Lind96,Lind98} constructed counterexamples to the local well-posedness for semilinear and quasilinear wave equations in three dimensions. Low-regularity local well-posedness for the quasilinear wave equation was proven by Tataru-Smith \cite{tataru}. They showed that for $n$ dimensional quasilinear wave equations of the form $\partial_t^2{\varphi}-\Delta {\varphi}={\partial{\varphi}}\partial^2{\varphi}$, the Cauchy problems are locally well-posed in $H^s(\mathbb{R}^n)$ with $s>n/2+7/4$ for $n=2$ and $s>(n+3)/2$ for $n=3,4,5$. See also Wang \cite{Wang}. Recently, Ohlmann \cite{Ohl} generalized Lindblad's result \cite{Lind98} to the 2D case and proved the ill-posedness of a 2D quasilinear wave equation in the logarithmic Sobolev space $H^{\frac{11}{4}}(\ln H)^{-\beta}$ with $\beta>\frac12$. This space is slightly more singular\footnote{There is a logarithmic loss between Ohlmann's result and the desired $H^{\frac{11}{4}}$ ill-posedness.} than $H^{\frac{11}{4}}$. In this paper, we derive the desired $H^{\frac{11}{4}}$ ill-posedness for 2D elastic wave system. Our approach and conclusion also apply to the 2D quasilinear wave equation. With no symmetry assumption, Alinhac \cite{Alinhac99,Alinhac99II} proved singularity formation for solutions to quasilinear wave equations in more than one spatial dimension via a Nash-Moser iteration scheme. This approach does not reveal information beyond the first blow-up point. In \cite{christodoulou10}, Christodoulou developed a geometric approach and provided a detailed understanding and a complete description of shock formation for 3D irrotational Euler's equations. This remarkable work was extended to a large class of equations, see \cite{Ross,speckbk,Speck16,Speck18,Speck-luk,Speck-luk2,miao,christodoulou-miao}.

Besides his well-known contributions in general relativity and in compressible fluids, Christodoulou also studied problems arising from mathematical physics. Our approach in this paper is based on a result of Christodoulou-Perez in \cite{christodoulou}, where they studied the propagation of electromagnetic waves in nonlinear crystals. Under planar symmetry, these electromagnetic waves satisfy a first-order genuinely nonlinear and strictly hyperbolic system. By revisiting and further extending John's work \cite{john74}, Christodoulou-Perez gave a more detailed description of the behaviour of solutions at the shock-formation time. In particular, they revealed that the blow-up of the first derivatives in John's work is equivalent to the vanishing of the inverse density for the characteristics. They also provided higher order estimates and a precise estimate for the blow-up time. Here we apply and extend this approach to the 2D elastic wave system and prove its desired low-regularity ill-posedness.

\subsection{Background and main results} \label{eq and main thm}
To describe the motion of an elastic body in 2D, we denote the time-dependent material deformation as $\Psi:\mathbb{R}^2\times\mathbb{R}\to\mathbb{R}^2$. Initially at $t=0$, $\Psi=(\Psi_1,\Psi_2)(x_1,x_2,t)$ satisfies $\Psi(x_1,x_2,0)=(x_1,x_2)$. The deformation gradient is defined as $F:=\nabla_x\Psi$  with  $F^i_{j}:=\partial \Psi_i/\partial x_j.$

An elastic body is called homogeneous hyperelastic if there exists a storage energy function $\hat{W}(F)$, such that its equation of motion can be derived via applying the least action principle to
\begin{align*}
  \mathcal{L}:=\int\int_{\mathbb{R}^2}\Big(\frac12|\partial_t\Psi|^2-\hat{W}(F)\Big)dx_1dx_2dt.
\end{align*}
The corresponding Euler-Lagrange equation takes the following form:
\begin{equation}\label{el}
  \frac{\partial^2\Psi^i}{\partial t^2 }-\frac{\partial}{\partial{x_l}}\frac{\partial \hat{W}(F)}{\partial F^i_{l}}=0.
\end{equation}

In this paper, we consider the homogeneous, hyperelastic, isotropic elastic body. For such an object, the storage energy function $\hat{W}(F)$ depends only on the principal invariants of $FF^T$. For notational simplification, we denote $\hat{W}(F)=W(k_1,k_2)$ with $k_1,k_2$ being the principal invariants of $C=FF^T-I$. And
 \begin{equation}
   \begin{split}
     k_1&=\mu_1+\mu_2=\text{tr} C,\\
     k_2&=\mu_1\mu_2=\frac12[(\text{tr} C)^2-\text{tr} C^2],
   \end{split}
 \end{equation}
where $\mu_1,\ \mu_2$ are the eigenvalues of $C$. Let $U=(U_1,U_2)^T:=(\Psi_1,\Psi_2)^T-(x_1,x_2)^T$ be the displacement and $G:=\nabla U=F-I$ be the displacement gradient. By a direct check, we have
\begin{equation} \label{k1k2}
   \begin{cases}
     k_1=\text{tr} G+\text{tr} GG^T,\\
     k_2=2(\text{tr} G)^2-\left[\text{tr} (G^2)+\text{tr}GG^T\right]+2\left[\text{tr}G\cdot \text{tr}(GG^T)-\text{tr}(G^2G^T)\right]+\frac12[(\text{tr}GG^T)^2-\text{tr}((GG^T)^2)].
   \end{cases}
 \end{equation}
Rewrite the storage energy function with the Taylor's expansion theorem, we have
\begin{equation}\label{taylor}
  W(k_1,k_2)=\gamma_0+\gamma_1k_1+\frac12\gamma_{11}k_1^2+\gamma_2 k_2+\frac16\gamma_{111}k_1^3+\gamma_{12}k_1k_2+\cdots
\end{equation}
Here, the constant coefficients $\gamma_0,\gamma_1,\gamma_{11},\cdots$ represent the values of the corresponding partial derivatives of $W$ at $k_i=0$ ($i=1,2$). We impose the condition that $\gamma_1=0$. This condition indicates that the reference configuration is a stress-free state. And $4(\gamma_{11}+\gamma_2)$, $-2\gamma_2$ are called the Lam\'e constants and we require them to be positive. The positivity of Lam\'e constants immediately gives $-2\gamma_2>0$ and $4(\gamma_{11}+\gamma_2)-2\gamma_2=2(2\gamma_{11}+\gamma_2)>0$. Noting that $G=\nabla U$ and invoking \eqref{k1k2}\eqref{taylor} in \eqref{el}, we arrive at the elastic wave equations
\begin{equation} \label{wave}
    \partial_t^2 U-c_2^2\Delta U-(c_1^2-c_2^2)\nabla(\nabla\cdot U)={\text {nonlinear terms of } (\nabla U,\nabla^2U)},
\end{equation}
where
$$c_1^2=4\gamma_{11},\quad c_2^2=-2\gamma_2\quad \text{and}\quad c_1^2>c_2^2>0.$$

Denote the quadratic nonlinear terms in \eqref{wave} to be $N(\nabla U, \nabla^2 U)$. Other nonlinear terms in \eqref{wave} take form of $(\nabla U)^\alpha\cdot\nabla U\cdot\nabla^2U$ with $\alpha\geq1$. These higher-order nonlinear terms are negligible compared with the quadratic ones in $N$ of the form $\nabla U\cdot\nabla^2 U$. It is because in this paper, under the evolution of the constructed initial data we will construct, we would prove the second derivatives of the displacement tend to be infinite instantaneously, while the first derivatives of the displacement would remain small. Without loss of generality, we keep only the quadratic nonlinear terms of \eqref{wave} in $N(\nabla U,\nabla^2 U)$ and hence derive the following quasilinear wave system with multiple wave-speeds:
\begin{equation}\label{y1}
    \partial_t^2 U-c_2^2\Delta U-(c_1^2-c_2^2)\nabla(\nabla\cdot U)=N(\nabla U,\nabla^2U)
\end{equation}
with $N(\nabla U,\nabla^2U)\sim\nabla U\nabla^2U$. Specifically, the nonlinearity takes the following form
\begin{equation}\label{elasqc}
\begin{split}
  &N(\nabla U,\nabla^2U)\\
  =&\sigma_0\nabla(\nabla\cdot U)^2+\sigma_1\big[\nabla(\nabla^\bot\cdot U)^2+2\nabla^\bot(\nabla\cdot U\nabla^\bot\cdot U)\big]\\
  &+\sigma_2\nabla Q_{12}(U_1,U_2)+\sigma_2\big(Q_{12}(\nabla\cdot U,U_2),Q_{12}(U_1,\nabla\cdot U)\big)^\top,
\end{split}
\end{equation}
where $\nabla^\bot=(-\partial_2,\partial_1)$,  $Q_{12}(f,g)=\partial_1f\partial_2g-\partial_2f\partial_1g$ and
\begin{equation}
\begin{array}{llll}
&\sigma_0=6\gamma_{11}+4\gamma_{111},\\
&\sigma_1=2(\gamma_{11}-\gamma_{12}),\\
&\sigma_2=2(\gamma_2-2\gamma_{11}+4\gamma_{12}).
\end{array}
\end{equation}

In this paper, we study elastic wave system \eqref{y1} and we deal with the general case
\begin{equation}\label{sigma}
  \sigma_0\sigma_1\neq 0.
\end{equation}
Note that the null conditions fail. Our main result is stated as below.
\begin{thm} \label{3D}
The Cauchy problem of the 2D elastic wave equations \eqref{y1} is ill-posed in $H^{\frac{11}{4}}(\mathbb{R}^2)$ in the following sense: we construct a family of compactly supported smooth initial data $(U_0^{(\eta)}, U_{t0}^{(\eta)})$ with
 $${\|U_0^{(\eta)}\|}_{{\dot{H}}^{\frac{11}{4}}(\mathbb{R}^2)}+{\|U_{t0}^{(\eta)}\|}_{{\dot{H}}^{\frac{7}{4}}(\mathbb{R}^2)} \to 0,\quad \text{as} \quad \eta\to0.$$
 Let $T_\eta^*$ be the largest time (a shock formation time) such that the Cauchy problem of \eqref{y1} $($with a general condition \eqref{sigma}$)$ admits a solution $U_\eta\in C^\infty(\mathbb{R}^2\times[0,T_\eta^*) )$. Then, as $\eta \to 0$, we have $T_\eta^* \to 0$.

Moreover, for each $\eta$, at the shock formation time $T^*_{\eta}$, the $H^2$ norm of the solution in a spatial region $\Omega_{T_\eta^*}$ blows up:
\begin{equation*}
  \|U_\eta(\cdot,T_\eta^*)\|_{H^2(\Omega_{T_\eta^*})}=+\infty.
\end{equation*}
\end{thm}
\begin{rmk}
In particular, by setting $u_2=0$ in \eqref{main} below, one can see that our $H^\frac{11}{4}$ ill-posedness also holds for the 2D quasilinear scalar wave equation:
\begin{equation} \label{single}
  \partial_t^2\varphi-\Delta \varphi=\partial_{x_1}(\partial_{x_1} \varphi)^2.
\end{equation}
Recall that Ohlmann \cite{Ohl} showed the ill-posedness of 2D quasilinear wave equation in the logarithmic Sobolev space $H^{\frac{11}{4}}(\ln H)^{-\beta}$ with $\beta>\frac12$. For \eqref{single}, our result here overcomes the logarithmic loss and reaches the desired regularity as suggested by Smith-Tataru \cite{tataru}.
\end{rmk}

\subsection{Strategy of the proof}
As in \cite{Lind96,Lind98} by Lindblad, we prove the ill-posedness under planar symmetry. Let $u(x,t)=(u_1(x,t),u_2(x,t))$ be a planar symmetric solution to \eqref{y1}. The elastic waves $(u_1,u_2)$ then obey
\begin{align} \label{main}
   \left\{\begin{array}{lll}
    \partial_t^2u_1-c_1^2\partial_x^2u_1=\sigma_0\partial_x(\partial_xu_1)^2+\sigma_1\partial_x(\partial_xu_2)^2,\\
    \partial_t^2u_2-c_2^2\partial_x^2u_2=2\sigma_1\partial_x(\partial_xu_1\partial_xu_2).
    \end{array}\right.
\end{align}
where $c_1>c_2>0$ and $\sigma_0,\sigma_1$ are constants. We then rewrite \eqref{main} as a $4\times 4$ first-order hyperbolic system:
\begin{equation}\label{1order}
        \partial_t\Phi+A(\Phi)\partial_x\Phi=0,
\end{equation}
where $\Phi=(\phi_1,\phi_2,\phi_3,\phi_4)^T=(\partial_xu_1,\partial_xu_2,\partial_tu_1,\partial_tu_2)^T$. Here we algebraically decompose the system via employing John's wave-decomposition method in \cite{john74}. We then carefully explore the structures of the corresponding system in a similar manner as Christodoulou and Perez did in \cite{christodoulou} for the nonlinear crystal optics.

In particular, based on the decomposition of waves, we trace the evolution of $\rho_i$ ($i=1,\cdots,4$), the so-called inverse densities of the $i^{\text{th}}$ characteristics. By constructing suitable initial data, we derive a positive uniform lower bound for $\{\rho_i\}_{i=2,3,4}$. While, we have that $\rho_1(z_0,t)\to 0$ as $t\to T_\eta^*$. Here, $T_\eta^*<+\infty$ is the first shock formation time.

We then calculate the $H^2(\mathbb{R}^2)$ norm of the solution to \eqref{y1} at time $T_\eta^*$. In a suitable constructed spatial region $\Omega_{T^*_\eta}$, we deduce that $\|U(\cdot,T^*_\eta)\|_{H^2(\Omega_{T^*_\eta})}^2=+\infty$. Furthermore, as $\eta \to 0$, we have $T_\eta^*\to 0$. This implies the desired ill-posedness result. We further unearth the hidden mechanism driving this ill-posedness. And we find that both the blow-up of $H^2$ norm and the vanishing of $T_\eta^*$ are driven by the shock formation: $\rho_1(z_0,t)\to 0$ as $t\to T_\eta^*.$

\subsection{Organization of the paper}
The paper is organized as follows. Section \ref{pre} is a preliminary part, where we rewrite the elastic wave system as a first order quasilinear hyperbolic system and employ decomposition of waves. In Section \ref{2d data}, we construct the $H^\frac{3}{4}$ initial data for the decomposed system, which corresponds to the $H^\frac{11}{4}$ initial data for the elastic wave system. Section \ref{aprior} is devoted to the proof of the shock formation. We apply Christodoulou-Perez's approach to our elastic wave system and depict the detailed solutions' behaviours up to the first shock formation time. In Section \ref{pfill}, we prove the $H^\frac{11}{4}$ ill-posedness and show that this ill-posedness is driven by the shock formation.

\section{Wave decomposition} \label{pre}
In this section, we employ an algebraic approach--wave decomposition to study the dynamics of the elastic wave system \eqref{y1}. Assuming that $U=(U_1,U_2)$ is a solution to \eqref{y1}, under plane symmetry, we have
\begin{equation*}
\begin{split}
U_1(x_1,x_2,t)=u_1(x_1,t),\ U_2(x_1,x_2,t)=u_2(x_1,t).
\end{split}
\end{equation*}
For notational simplicity, we denote $x_1$ by $x$ in below. Set $u(x,t)=(u_1(x,t),u_2(x,t))$ to be a planar symmetric solution to \eqref{y1}, then it satisfies a quasilinear wave system with multiple speeds:
\begin{align}\label{2d1}
    \left\{\begin{array}{lll}
    \partial_t^2u_1-c_1^2\partial_x^2u_1=\sigma_0\partial_x(\partial_xu_1)^2+\sigma_1\partial_x(\partial_xu_2)^2,\\
    \partial_t^2u_2-c_2^2\partial_x^2u_2=2\sigma_1\partial_x(\partial_xu_1\partial_xu_2).
    \end{array}\right.
\end{align}
Now define
\begin{align*}
  \phi_1:=\partial_xu_1,\quad \phi_2:=\partial_xu_2,\quad \phi_3:=\partial_tu_1,\quad \phi_4:=\partial_tu_2.
\end{align*}
and let $\Phi:=(\phi_1,\phi_2,\phi_3,\phi_4)^T$. We hence rewrite system \eqref{2d1} as
\begin{align}\label{2d2}
  \partial_t\Phi+A(\Phi)\partial_x\Phi=0,
\end{align}
where
\begin{align*}
  \begin{split}
    A(\Phi)=\left(
              \begin{array}{cccccc}
                0 & 0  &-1 & 0  \\
                0 & 0 &0  & -1  \\
                -(c_1^2+2\sigma_0\phi_1) & -2\sigma_1\phi_2 &0  & 0  \\
                -2\sigma_1\phi_2 & -(c_2^2+2\sigma_1\phi_1)&0  & 0  \\
              \end{array}
            \right).
  \end{split}
\end{align*}

{\it  Eigenvalues and Eigenvectors of $A(\Phi)$.} We first compute $A(\Phi)$'s eigenvalues and choose the corresponding eigenvectors. Via calculations, the eigenvalues of $(A(\Phi))_{4\times4}$ are
\begin{align}\label{118}
  \begin{array}{lll}
    \lambda_1=&\sqrt{\frac12(a+b)+\frac12\sqrt{(a-b)^2+4c^2}},\\
    \lambda_2=&\sqrt{\frac12(a+b)-\frac12\sqrt{(a-b)^2+4c^2}},\\
    \lambda_3=&-\sqrt{\frac12(a+b)-\frac12\sqrt{(a-b)^2+4c^2}},\\
    \lambda_4=&-\sqrt{\frac12(a+b)+\frac12\sqrt{(a-b)^2+4c^2}},
  \end{array}
\end{align}
with
\begin{equation}\label{abc}
 a=c_1^2+2\sigma_0\phi_1,\quad b=c_2^2+2\sigma_1\phi_1,\quad c=2\sigma_1\phi_2.
\end{equation}
A quasilinear hyperbolic system is called strictly hyperbolic if the eigenvalues of its coefficient matrix are all distinct. Note that the elastic wave system \eqref{2d2} is strictly hyperbolic for $|\Phi|<2\kappa$ with $\kappa$ small enough. In particular, for $|\Phi|<2\kappa$, it holds that
\begin{align}\label{ei}
    \lambda_4(\Phi)<\lambda_3(\Phi)<\lambda_2(\Phi)<\lambda_1(\Phi).
\end{align}
We then design the right eigenvectors as:
\begin{align}\label{right}
\begin{array}{lll}
 & r_1=\left(
                \begin{array}{c}
                  \frac{\lambda_1^2-b}{2\sigma_1} \\
                  \phi_2 \\
                  -\frac{\lambda_1(\lambda_1^2-b)}{2\sigma_1} \\
                  -\lambda_1\phi_2 \\
                \end{array}
              \right),\ r_2=\left(
                                     \begin{array}{c}
                                       \frac{\lambda_2^2-b}{2\sigma_1} \\
                                       \phi_2 \\
                                       -\frac{\lambda_2(\lambda_2^2-b)}{2\sigma_1} \\
                                       -\lambda_2\phi_2\\
                                     \end{array}
                                   \right),\ r_3=\left(
                                                   \begin{array}{c}
                                                     \frac{\lambda_2^2-b}{2\sigma_1} \\
                                                     \phi_2 \\
                                                     \frac{\lambda_2(\lambda_2^2-b)}{2\sigma_1} \\
                                                     \lambda_2\phi_2 \\
                                                   \end{array}
                                                 \right),\ r_4=\left(
                                                                 \begin{array}{c}
                                                                    \frac{\lambda_1^2-b}{2\sigma_1} \\
                                                                   \phi_2 \\
                                                                   \frac{\lambda_1(\lambda_1^2-b)}{2\sigma_1}\\
                                                                   \lambda_1\phi_2  \\
                                                                 \end{array}
                                                               \right).
\end{array}
\end{align}
We set the corresponding left eigenvectors to be the dual of the right ones, i.e.,
\begin{equation}\label{lr}
  l_i(\Phi)\cdot r_j(\Phi)=\delta_{ij}.
\end{equation}
And they are
\begin{align}\label{left}
\begin{array}{lll}
  l_1=&\frac1{K}\Big(\frac{\lambda_1^2-b}{2\sigma_1},\phi_2,-\frac{\lambda_1^2-b}{2\sigma_1\lambda_{1}},-\frac{\phi_2}{\lambda_{1}}\Big),\
  &l_2=\frac1{N}\Big(\frac{\lambda_2^2-b}{2\sigma_1},\phi_2,-\frac{\lambda_2^2-b}{2\sigma_1\lambda_{2}},-\frac{\phi_2}{\lambda_{2}}\Big),\\ l_3=&\frac1{N}\Big(\frac{\lambda_2^2-b}{2\sigma_1},\phi_2,\frac{\lambda_2^2-b}{2\sigma_1\lambda_{2}},\frac{\phi_2}{\lambda_{2}}\Big),\ &l_4=\frac1{K}\Big(\frac{\lambda_1^2-b}{2\sigma_1},\phi_2,\frac{\lambda_1^2-b}{2\sigma_1\lambda_{1}},\frac{\phi_2}{\lambda_{1}}\Big),
\end{array}
\end{align}
where
 \begin{align*}
   K=&\frac{(a-b)^2+4c^2}{4\sigma_1^2}+\frac{(a-b)\sqrt{(a-b)^2+4c^2}}{4\sigma_1^2},\\
   N=&\frac{(a-b)^2+4c^2}{4\sigma_1^2}-\frac{(a-b)\sqrt{(a-b)^2+4c^2}}{4\sigma_1^2}=\frac{(\lambda_2^2-b)^2+c^2}{2}.
 \end{align*}
Note that $\lambda_2^2-b=0$ if and only if $\phi_2=0$. One can also verify that the derivatives with respect to $\phi_2$ satisfy $\partial_{\phi_2}(\lambda_2^2-b)=O(\phi_2)$ and $\partial^2_{\phi_2\phi_2}(\lambda_2^2-b)=O(1)$. Hence, the leading order of $\lambda_2^2-b$ is $O(\phi_2^2)$.

{\it Characteristic and bi-characteristic coordinates.}  As in \cite{an, an2,christodoulou}, we employ the characteristic coordinates and bi-characteristic coordinates for \eqref{2d2}.  The  {\bf\textit{ $i^{\text{th}}$ characteristic $\mathcal{C}_i(z_i)=\Big\{\Big(X_i(z_i,t),t\Big): 0\leq t\leq T\Big\}$}} is defined via the flow map, which originates from $z_i$ at $t=0$ and satisfies:
\begin{align}\label{flow}
  \left\{\begin{array}{ll}
  \frac{\partial}{\partial t}X_i(z_i,t)=\lambda_i\big(\Phi(X_i(z_i,t),t)\big),\quad t\in[0,T],\\
  X_i(z_i,0)=z_i.
  \end{array}\right.
\end{align}
For any $(x,t)\in \mathbb{R}\times [0,T]$, there is a unique $(z_i,s_i)\in \mathbb{R}\times [0,T]$ such that $(x,t)=\big(X_i(z_i,s_i),s_i\big)$ with $X_i$ satisfying \eqref{flow}. We define this $(z_i,s_i)$ to be {\bf\textit{the characteristic coordinate}} of $(x,t)$.
We then introduce the bi-characteristic coordinates to describe the intersection of two transversal characteristics $\mathcal{C}_i(z_i)$ and $\mathcal{C}_j(z_j)$ when $i\neq j$. For any given $(x,t)\in \mathbb{R}\times [0,T]$, there is a unique pair of $(z_i,z_j)$ such that the characteristics $\mathcal{C}_i$ and $\mathcal{C}_j$ meet at $(x,t)$. To distinguish with the characteristic coordinates above, we denote $(z_i,z_j)$ by $(y_i,y_j)$, and define $(y_i,y_j)$ as {\bf\textit{the bi-characteristic coordinate}}. In particular, these coordinates obey
\begin{equation}\label{bichar}
  (x,t)=\big(X_i(y_i,t'(y_i,y_j)),t'(y_i,y_j)\big)=\big(X_j(y_j,t'(y_i,y_j)),t'(y_i,y_j)\big).
\end{equation}
To describe the compression among the $i^{\text{th}}$ characteristics, following \cite{christodoulou}, we employ a geometric quantity $\rho_i$ (the inverse density of the $i^{\text{th}}$ characteristics)
\begin{equation}\label{dense}
\rho_i:=\partial_{z_i}X_i.
\end{equation}
Following from \eqref{flow}, it holds
\begin{equation}\label{319}
  \rho_i(z_i,0)=1.
\end{equation}
Via direct calculations, we have that the coordinate transformations comply with the following rules:
\begin{equation}
  \partial_{z_i}=\rho_i\partial_x,\quad \partial_{s_i}=\lambda_i\partial_x+\partial_t,
\quad
  \partial_{y_i}t'=\frac{\rho_i}{\lambda_j-\lambda_i},\quad\partial_{y_j}t'=\frac{\rho_j}{\lambda_i-\lambda_j},
\end{equation}
\begin{equation}\label{biytoz}
  \partial_{y_i}=\frac{\rho_i}{\lambda_j-\lambda_i}\partial_{s_j}=\partial_{z_i}+\frac{\rho_i}{\lambda_j-\lambda_i}\partial_{s_i},
\end{equation}
and
\begin{equation}
  dx=\rho_idz_i+\lambda_ids_i,\quad dt=ds_i,\quad   dz_i=dy_i,\quad dz_j=dy_j,
\end{equation}
\begin{equation}
  dx=\frac{\rho_i\lambda_j}{\lambda_j-\lambda_i}dy_i+\frac{\rho_j\lambda_i}{\lambda_i-\lambda_j}dy_j,\quad dt=\frac{\rho_i}{\lambda_j-\lambda_i}dy_i+\frac{\rho_j}{\lambda_i-\lambda_j}dy_j.
\end{equation}

{\it Wave decomposition.}  Now we are ready to decompose system \eqref{2d2} via the decomposition of waves. This approach was first introduced by John \cite{john74} for genuinely nonlinear strictly hyperbolic system. Christodoulou and Perez \cite{christodoulou} revisited this approach with a geometric view and sharpened John's result by exhibiting the details of the shock formation. See \cite{an,an2,liu} for more applications of this method. For a fixed $i \in \{1,2,3,4\}$, define
\begin{equation}\label{wi}
  w_i:=l_i\partial_x\Phi
  \end{equation}
and
  \begin{equation}
    v_i:=\rho_iw_i.
\end{equation}
Then, by \eqref{lr}, we have
\begin{equation}\label{phix}
  \partial_x\Phi=\sum_{k=1}^4 w_kr_k.
\end{equation}
The above formula serves as a \underline{decomposition} of $\partial_x \Phi$. One can verify that $(\rho_i,w_i,v_i)$ satisfy the following transport equations
\begin{align}
  \partial_{s_i}\rho_i=&c_{ii}^iv_i+\Big(\sum_{m\neq i}c_{im}^iw_m\Big)\rho_i,\label{eqrho}\\
  \partial_{s_i}w_i=&-c_{ii}^iw_i^2+\Big(\sum_{m\neq i}(-c_{im}^i+\gamma_{im}^i)w_m\Big)w_i+\sum_{m\neq i,k\neq i\atop m\neq k}\gamma_{km}^iw_kw_m,\label{eqw}\\
  \partial_{s_i}v_i=&\Big(\sum_{m\neq i}\gamma_{im}^iw_m\Big)v_i+\sum_{m\neq i,k\neq i\atop m\neq k}\gamma_{km}^iw_kw_m\rho_i,\label{eqv}
\end{align}
with $\partial_{s_i}=\lambda_i\partial_x+\partial_t$ and
\begin{align}
&c_{im}^i=\nabla_\Phi\lambda_i\cdot r_m,\label{coec}\\
  &\gamma_{im}^i=-(\lambda_i-\lambda_m)l_i \cdot(\nabla_\Phi r_i \cdot r_m-\nabla_\Phi r_m \cdot r_i),\quad m\neq i,\label{coeg1}\\
  &\gamma_{km}^i=-(\lambda_k-\lambda_m)l_i \cdot (\nabla_\Phi r_k \cdot r_m), \qquad\qquad\qquad k\neq i,\  m\neq i.\label{coeg2}
\end{align}

{\it Bounds for the coefficients. } We estimate the coefficients $c_{im}^i, \gamma_{im}^i, \gamma_{km}^i $ in the equations \eqref{eqrho}-\eqref{eqv}. Assume $\Phi\in B_{2\kappa}^4(0)$ for some fixed positive constant $\kappa\ll1$. With the eigenvectors chosen in \eqref{right} and \eqref{left}, we have that these coefficients are uniformly bounded. By \eqref{118} and definitions of $a,b,c$ in \eqref{abc}, for $\Phi\in B_{2\kappa}^4(0)$, we have
\begin{equation*}
  a-b\sim c_1^2-c_2^2,\quad \lambda_1^2-b\sim c_1^2-c_2^2,\quad c\sim 0.
\end{equation*}
From the definitions of these coefficients, the only potential singularity would come from the factor $N$ in the left eigenvectors $l_2$ and $l_3$ in \eqref{left}. These left eigenvectors are used to define $\gamma^2_{km}$ and $\gamma^3_{km}$ in \eqref{coeg1} and \eqref{coeg2}. One can further verify that all the (potentially singular) coefficients $\gamma^2_{km}$, $\gamma^3_{km}$ are of $O(1)$. For example, we have
\begin{eqnarray*}
  \gamma^2_{13} &=& -\frac{\lambda_1-\lambda_3}{N}\left[\frac{(\lambda_2^2-b)^2\partial_{\phi_1}(\lambda_1^2-b)}{8\sigma_1^3}+\frac{(\lambda_2^2-b)\phi_2\partial_{\phi_2}(\lambda_1^2-b)}{4\sigma_1^2}+\phi_2^2\right.\\
   &&+\frac{\lambda_1(\lambda_2^2-b)^2\left[\partial_{\phi_1}(\lambda_1^2-b)+(\lambda_1^2-b)\partial_{\phi_1}\lambda_1\right]}{8\sigma_1^3\lambda_2}+\frac{(\lambda_2^2-b)\left[\lambda_1\partial_{\phi_2}(\lambda_1^2-b)+(\lambda_1^2-b)\partial_{\phi_2}\lambda_1\right]}{4\sigma_1^2\lambda_2}\\
   &&\left.+\frac{(\lambda_2^2-b)\phi_2^2\partial_{\phi_1}\lambda_2}{2\sigma_1\lambda_2}+\phi_2^3\left(1+\frac{\partial_{\phi_2}\lambda_2}{\lambda_2}\right)\right]\\
   &=&\frac{1}{N}\left(O((\lambda_2^2-b)^2)+O(\phi_2^2)\right)=O(1).
\end{eqnarray*}
The last equality holds because of $N=\frac{(\lambda_2^2-b)^2+c^2}{2}$ and $c^2=4\sigma_1^2 \phi_2^2$. Other potentially singular coefficients in $\gamma^2_{km}$, $\gamma^3_{km}$ can be estimated similarly. They are all of order $O(1)$.

We also note that the $1^{\text{st}}$ and $4^{\text{th}}$ characteristics are {\bf\textit{genuinely non-linear}} in the sense of Lax, i.e., $c^1_{11},c^4_{44}\neq0$. Denote
\begin{equation*}
  \Delta=(a-b)^2+4c^2.
\end{equation*}
Then, it holds that
\begin{align*}
  &\nabla_\Phi\lambda_1=-\nabla_\Phi\lambda_4=\Big(\frac{(\sigma_0+\sigma_1)\sqrt{\Delta}+(a-b)(\sigma_0-\sigma_1)}{2\lambda_1\sqrt{\Delta}},\frac{4\sigma_1^2\phi_2}{\lambda_1\sqrt{\Delta}},0,0\Big),\\
  &\nabla_\Phi\lambda_2=-\nabla_\Phi\lambda_3=\Big(\frac{(\sigma_0+\sigma_1)\sqrt{\Delta}-(a-b)(\sigma_0-\sigma_1)}{2\lambda_2\sqrt{\Delta}},-\frac{4\sigma_1^2\phi_2}{\lambda_2\sqrt{\Delta}},0,0\Big).
\end{align*}
Therefore, we have
\begin{equation*}
  c_{11}^1(\Phi)=-c_{44}^4(\Phi)=\nabla_\Phi\lambda_1 \cdot r_1=\frac{2\sigma_0(a-b)(\lambda_1^2-b)+(2\sigma_0+6\sigma_1)c^2}{4\sigma_1\lambda_1\sqrt{\Delta}},
\end{equation*}
\begin{equation*}
  c_{22}^2(\Phi)=-c_{33}^3(\Phi)=\nabla_\Phi\lambda_2 \cdot r_2=-\frac{2\sigma_0(a-b)(\lambda_2^2-b)-(2\sigma_0+6\sigma_1)c^2}{4\sigma_1\lambda_2\sqrt{\Delta}}.
\end{equation*}
Notice that
\begin{equation*}
  a-b\sim c_1^2-c_2^2,\quad \lambda_1^2-b\sim c_1^2-c_2^2,\quad c\sim 0.
\end{equation*}
This implies that $c_{22}^2, c_{33}^3$ are small perturbations around zero. And $c_{11}^1, c_{44}^4$ are small perturbations around $\frac{\sigma_0(c_1^2-c_2^2)}{2\sigma_1c_1},-\frac{\sigma_0(c_1^2-c_2^2)}{2\sigma_1c_1}$, which are both away from zero. In the following, without loss of generality, we consider the elastic body satisfying
\begin{equation}
c^1_{11}(0)=\frac{\sigma_0(c_1^2-c_2^2)}{2\sigma_1c_1}<0.
\end{equation}

\section{Construction of initial data} \label{2d data}
In this section, we construct $H^\frac34$ initial data for the Cauchy problem of the decomposed system \eqref{eqrho}-\eqref{eqv}. These initial data then yield $H^\frac{11}{4}$ initial data for the 2D elastic wave system \eqref{y1}. For any fixed $x_2$, define
\begin{equation}\label{dataw0}
\hat{w}_1^{(\eta)}(x,x_2)=w_1^{(\eta)}(x,x_2,0)=\theta |\ln (x)|^\alpha \mathcal{X}_\eta(x)\psi\Big(\frac{|\ln(x)|^\delta x_2}{\sqrt{x}}\Big)
\end{equation}
with $0<\theta\ll 1$, $0<\alpha<\frac12$ and $\delta>0$ being three constants. Here, $\mathcal{X}_\eta(x)=\mathcal{X}(\frac{x}{\eta})$ and we require
$$\mathcal{X}(x)=\left\{\begin{split}
&1,\quad x\in[\frac65,\frac95],\\
&0,\quad x\in(-\infty, 1]\cup [2,+\infty),
\end{split}\right.
\qquad \psi(x)=\left\{\begin{split}
&1,\quad |x|\leq \frac14,\\
&0,\quad |x|\geq\frac12.
\end{split}\right.
$$
Moreover, we set
\begin{equation} \label{data2}
\hat{w}_k^{(\eta)}(x,x_2)=w_k^{(\eta)}(x,x_2,0)=\theta^2 \mathcal{X}_\eta(x)\psi\Big(\frac{|\ln(x)|^\delta x_2}{\sqrt{x}}\Big)\quad\text{for}\quad k=2,3,4.
\end{equation}
For the above constructed initial data, we have
\begin{equation}\label{W0}
  W_0^{(\eta)}:=\max_i\sup_{(z_i,x_2)}|w_i^{(\eta)}(z_i,x_2,0)|=w_1^{(\eta)}(z_0,x_2,0)>0.
\end{equation}
The construction of the above initial data here is inspired by Lindblad \cite{Lind96,Lind98}, Ohlmann \cite{Ohl} and An-Chen-Yin \cite{an}. Lindblad picked the function $|\ln (x)|^\alpha$ in \cite{Lind96} to generate low-regularity initial data that yield ill-posedness for 3D wave equations in \cite{Lind96,Lind98}. Ohlmann \cite{Ohl} generalized Lindblad's result to 2D case by properly designing a test function $\psi$ which determines the geometry of the domain of future dependence. We modify Ohlmann's data by including a cutoff function $\chi$ as designed in \cite{an}.

With the above constructed initial data, one can verify
\begin{lem}
For $\hat{w}_1^{(\eta)}(x,x_2)$ defined in \eqref{dataw0}. If $2\alpha-\delta<0$, then we have $\hat{w}_1^{(\eta)}\in \dot{H}^{\frac34}(\mathbb{R}^2)$ and $\|\hat{w}_1^{(\eta)}\|_{\dot{H}^{\frac34}} \to 0$ as $\eta \to 0$.
\end{lem}
{\it Proof.} A detailed proof for a more generalized construction is provided in \cite{an3} by us. Here we outline some key ideas for the proof of this lemma. We first calculate the Fourier transform of $\hat{w}_1^{(\eta)}$:
\begin{equation}\label{fw}
\begin{split}
\mathcal{F}(\hat{w}_1^{(\eta)})(\xi_1,\xi_2)
&=\int_{x}e^{-2\pi i x\xi_1} \theta |\ln (x)|^\alpha \mathcal{X}_\eta(x)dx\int_{x_2}e^{-2\pi i x_2\xi_2}\psi\Big(\frac{|\ln(x)|^\delta x_2}{\sqrt x}\Big)dx_2\\
&=\int_{\eta}^{2\eta}e^{-2\pi i x\xi_1} \theta |\ln (x)|^\alpha \mathcal{X}_\eta(x)\frac{\sqrt x}{|\ln x|^{\delta }}dx\int_{|y|\leq\frac12}e^{-2\pi i y\frac{\sqrt x}{|\ln x|^\delta }\xi_2}\psi(y)dy.
\end{split}
\end{equation}
To estimate the $\dot{H}^{\frac34}(\mathbb{R}^2)$ norm of the data, we divide the domain of integration into four regions:
\begin{equation*}
\begin{split}
& D_1=\{(\xi_1,\xi_2):|\xi_1|\leq\frac1\eta, |\xi_2|\leq\frac{|\ln\eta|^\delta}{\sqrt \eta}\},\ D_2=\{(\xi_1,\xi_2):|\xi_1|>\frac1\eta, |\xi_2|\leq\frac{|\ln\eta|^\delta}{\sqrt \eta}\},\\
& D_3=\{(\xi_1,\xi_2):|\xi_1|\leq\frac1\eta, |\xi_2|>\frac{|\ln\eta|^\delta}{\sqrt \eta}\},\  D_4=\{(\xi_1,\xi_2):|\xi_1|>\frac1\eta, |\xi_2|>\frac{|\ln\eta|^\delta}{\sqrt \eta}\}.\\
\end{split}
\end{equation*}

In domain $D_1$, from \eqref{fw}, we directly obtain
\begin{equation}\label{D1}
\begin{split}
\iint_{D_1}|\xi|^{\frac32}[\mathcal{F}(\hat{w}_1^{(\eta)})(\xi_1,\xi_2)]^2d\xi_1d\xi_2\lesssim\eta^{-\frac32}|\ln\eta|^\delta\eta^{-\frac32}\theta^2\eta^3|\ln\eta|^{2\alpha-2\delta}=\theta^2|\ln\eta|^{2\alpha-\delta}.
\end{split}
\end{equation}

Within the non-compact domain $D_2$, we have
\begin{equation}\label{F1}
\begin{split}
\mathcal{F}(\hat{w}_1^{(\eta)})(\xi_1,\xi_2)
=&\frac{C_0}{\xi_1}\int_{\eta}^{2\eta}\theta |\ln (x)|^\alpha \mathcal{X}_\eta(x)\frac{\sqrt x}{|\ln x|^{\delta }}d(e^{-2\pi i x\xi_1} )\int_{|y|\leq\frac12}e^{-2\pi i y\frac{\sqrt x}{|\ln x|^\delta }\xi_2}\psi(y)dy\\
=&\frac{C_1}{\xi_1}\underbrace{\int_{\eta}^{2\eta}e^{-2\pi i x\xi_1} \theta \frac{|\ln (x)|^{\alpha-1}}{x} \mathcal{X}_\eta(x)\frac{\sqrt x}{|\ln x|^{\delta }}dx\int_{|y|\leq\frac12}e^{-2\pi i y\frac{\sqrt x}{|\ln x|^\delta }\xi_2}\psi(y)dy}_{I_1}\\
&+\frac{C_2}{\xi_1}\underbrace{\int_{\eta}^{2\eta}e^{-2\pi i x\xi_1} \theta |\ln (x)|^\alpha \mathcal{X}'_\eta(x)\frac{\sqrt x}{|\ln x|^{\delta }}dx\int_{|y|\leq\frac12}e^{-2\pi i y\frac{\sqrt x}{|\ln x|^\delta }\xi_2}\psi(y)dy}_{I_2}\\
&+\frac{C_3}{\xi_1}\underbrace{\int_{\eta}^{2\eta}e^{-2\pi i x\xi_1} \theta |\ln (x)|^\alpha \mathcal{X}_\eta(x)\frac{1}{\sqrt x|\ln x|^{\delta }}dx\int_{|y|\leq\frac12}e^{-2\pi i y\frac{\sqrt x}{|\ln x|^\delta }\xi_2}\psi(y)dy}_{I_3}\\
&+\frac{C_4}{\xi_1}\underbrace{\int_{\eta}^{2\eta}e^{-2\pi i x\xi_1} \theta |\ln (x)|^\alpha \mathcal{X}_\eta(x)\frac{\sqrt x}{x|\ln x|^{\delta +1}}dx\int_{|y|\leq\frac12}e^{-2\pi i y\frac{\sqrt x}{|\ln x|^\delta }\xi_2}\psi(y)dy}_{I_4}\\
&+\frac{C_5}{\xi_1}\underbrace{\int_{\eta}^{2\eta}e^{-2\pi i x\xi_1} \theta |\ln (x)|^\alpha \mathcal{X}_\eta(x)\frac{1}{|\ln x|^{2\delta }}dx\int_{|y|\leq\frac12}e^{-2\pi i y\frac{\sqrt x}{|\ln x|^\delta }\xi_2}\xi_2y\psi(y)dy}_{I_5}\\
&+\frac{C_6}{\xi_1}\underbrace{\int_{\eta}^{2\eta}e^{-2\pi i x\xi_1} \theta |\ln (x)|^\alpha \mathcal{X}_\eta(x)\frac{1}{|\ln x|^{2\delta+1 }}dx\int_{|y|\leq\frac12}e^{-2\pi i y\frac{\sqrt x}{|\ln x|^\delta }\xi_2}\xi_2y\psi(y)dy}_{I_6},
\end{split}
\end{equation}
where $C_i$ are uniform constants with $i=0,\cdots,6$. Because
\begin{align*}
\begin{array}{lll}
&|I_1|\lesssim\theta\sqrt\eta|\ln\eta|^{\alpha-\delta},\qquad&|I_2|\lesssim\theta\sqrt\eta|\ln\eta|^{\alpha-\delta},\\
&|I_3|\lesssim\theta\sqrt\eta|\ln\eta|^{\alpha-\delta},\qquad&|I_4|\lesssim\theta\sqrt\eta|\ln\eta|^{\alpha-\delta-1},\\
&|I_5|\lesssim\theta\eta|\ln\eta|^{\alpha-2\delta},\qquad&|I_6|\lesssim\theta\eta|\ln\eta|^{\alpha-2\delta-1},
\end{array}\end{align*}
it then follows that
\begin{equation}
|\mathcal{F}(\hat{w}_1^{(\eta)})(\xi_1,\xi_2)|\lesssim \frac{1}{|\xi_1|}\theta\eta^{\frac12}|\ln\eta|^{\alpha-\delta}.
\end{equation}
Integrating by parts twice with respect to $x$,  we derive
\begin{equation}\label{intep}
|\mathcal{F}(\hat{w}_1^{(\eta)})(\xi_1,\xi_2)|\lesssim \frac{1}{|\xi_1|^2}\theta\eta^{-\frac12}|\ln\eta|^{\alpha-\delta}.
\end{equation}
Thus, in $D_2$, it holds
\begin{equation}\label{D2}
\begin{split}
\iint_{D_2}|\xi|^{\frac32}[\mathcal{F}(\hat{w}_1^{(\eta)})(\xi_1,\xi_2)]^2d\xi_1d\xi_2&\lesssim\iint_{D_2}|\xi_1|^{\frac32}\theta^2\frac{1}{\xi_1^{4}}\eta^{-1}|\ln\eta|^{2\alpha-2\delta}d\xi_1d\xi_2\\
&=\theta^2\eta^{-1}|\ln\eta|^{2\alpha-2\delta}\frac{|\ln\eta|^\delta}{\sqrt\eta}\int_{|\xi_1|>\eta^{-1}}\xi_1^{-\frac52}d\xi_1\\
&\lesssim\theta^2|\ln\eta|^{2\alpha-\delta}.
\end{split}
\end{equation}
Similarly, we also have
\begin{equation}\label{D3}
\iint_{D_3}|\xi|^{\frac32}[\mathcal{F}(\hat{w}_1^{(\eta)})(\xi_1,\xi_2)]^2d\xi_1d\xi_2 \lesssim\theta^2|\ln\eta|^{2\alpha-\delta}
\end{equation}
and
\begin{equation}\label{D4}
\iint_{D_4}|\xi|^{\frac32}[\mathcal{F}(\hat{w}_1^{(\eta)})(\xi_1,\xi_2)]^2d\xi_1d\xi_2\lesssim\theta^2|\ln\eta|^{2\alpha-\delta}.
\end{equation}

In summary, from \eqref{D1} and \eqref{D2}-\eqref{D4}, we derive
\begin{equation}
\begin{split}
\|\hat{w}_1^{(\eta)}\|^2_{\dot{H}^{\frac34}(\mathbb{R}^2)}\lesssim \theta^2|\ln\eta|^{2\alpha-\delta}.
\end{split}
\end{equation}
When $2\alpha-\delta<0$, the right hand side of the above inequality is bounded. This completes the proof of this lemma.\hfill$\Box$

Next, by the formula of decomposition of waves \eqref{phix}, we have
\begin{equation} \label{combo2}
\partial_x\Phi(x,x_2,0)=\sum_{k=1}^4 {w}_i^{(\eta)}(x,x_2,0)r_i\big(\Phi(x,x_2,0)\big).
\end{equation}
Noting that our $r_i$ is Lipschitz continuous in $\Phi$, by a standard ODE argument, one can solve the above ODE system \eqref{combo2} for the initial data $\Phi(x,x_2,0)$. Furthermore, since $r_i(\Phi)\in L^\infty( B_{2\kappa}^4(0))$ and $w_1(x,x_2,0)\in H^{\frac34}(\mathbb{R}^2)$, there holds $\Phi(x,x_2,0)\in H^{\frac{7}{4}}(\mathbb{R}^2)$. Recall that $\Phi=(\partial_x U_1, \partial_x U_2, \partial_t U_1, \partial_t U_2)$. This implies that the initial data for the 2D elastic waves \eqref{y1} belong to $H^{\frac{11}{4}}(\mathbb{R}^2)$.

\section{Shock formation}\label{aprior}
In this section, we prove shock formation for the 2D elastic wave system. In particular, we shall derive estimates for the following norms:
\begin{align}
  S(t):=&\max_{i}\sup_{(z'_i,s'_i)\atop z'_i\in[\eta,2\eta],\ 0\leq s'_i\leq t}\rho_i(z'_i,s'_i),\\
   J(t):=&\max_{i}\sup_{(z'_i,s'_i)\atop z'_i\in[\eta,2\eta]\ 0\leq s'_i\leq t}|v_i(z'_i,s'_i)|,\\
  V(t):=&\max_i\sup_{(x',t')\notin\mathcal{R}_i,\atop 0\leq t'\leq t}|w_i(x',t')|,\\
  \bar{U}(t):=&\sup_{(x',t')\atop 0\leq t'\leq t}|\Phi(x',t')|.
\end{align}

Before proceeding to these estimates, we analyze some properties of the characteristics. We consider the case that the initial data are supported on $I_0=[\eta,2\eta]$. And we define the {\bf \textit{ $i^{\text{th}}$ characteristic strips $\mathcal{R}_i$}}:
\begin{equation}
\mathcal{R}_i=\cup_{z_i\in I_0}\mathcal{C}_i(z_i),\qquad i=1,2,3,4.
\end{equation}
These four characteristic strips $\mathcal{R}_1,\mathcal{R}_2,\mathcal{R}_3,\mathcal{R}_4$ will be totally separated after some time. Specifically, we first take the supremum and infimum of the eigenvalues
\begin{align*}
  \bar{\lambda}_i:=\sup_{\Phi\in B_{2\kappa}^4 (0)}\lambda_i(\Phi),\quad \underline{\lambda}_i:=\inf_{\Phi\in B_{2\kappa}^4(0)}\lambda_i(\Phi),\quad \text{for}\quad i=1,2,3,4
\end{align*}
and define
\begin{equation*}
  \sigma:=\min_{i<j\atop i,j\in\{1,2,,3,4\}}(\underline{\lambda}_i-\bar{\lambda}_j).
\end{equation*}
When $\kappa$ is sufficiently small, we have that $\sigma$ has a uniform positive lower bound. For $i\in\{1,2,3,4\}$, $z\in I_0$, it follows from \eqref{flow} that
\begin{equation*}
  z+\underline{\lambda}_i t\leq X_i(z,t)\leq z+\bar{\lambda}_i t.
\end{equation*}
Moreover, for all $i<j$, with $i,j\in\{1,2,3,4\}$, there holds
\begin{equation*}
   X_i(\eta,t)-X_j(2\eta,t) \geq (\eta+\underline{\lambda}_i t)-(2\eta+\bar{\lambda}_j t)=-\eta+(\underline{\lambda}_i-\bar{\lambda}_j)t\geq-\eta+\sigma t.
\end{equation*}
Note that the above difference is strictly positive when
\begin{equation}\label{t0}
  t>t_0^{(\eta)}:=\frac{\eta}\sigma>0.
\end{equation}
This implies that the four characteristic strips are well separated after $t_0^{(\eta)}$.
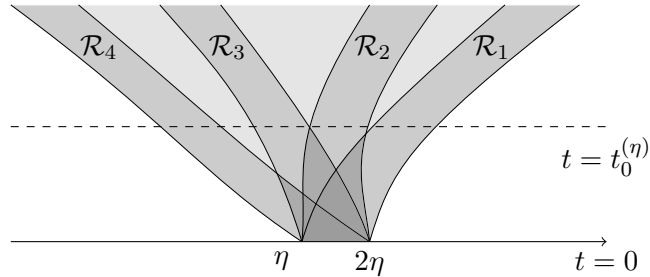
\begin{figure}[htb]
\centering
\begin{tikzpicture}[scale=0.9, fill opacity=0.5, draw opacity=1, text opacity=1]
\node [below]at(3.5,0){$2\eta$};
\node [below]at(2.2,0){$\eta$};

\filldraw[white, fill=gray!40] (3.5,0)..controls (2,1) and (1,2)..(-0.8,3.5)--(0.4,3.5)..controls (1,2.8) and (1.9,2)..(2.5,0);
\filldraw[white, fill=gray!40](3.5,0)..controls (3.2,1) and (2,2.6)..(1.3,3.5)--(3.5,3.5)..controls (2.2,1.5) and (2.6,1)..(2.5,0);
\filldraw[white, fill=gray!40](3.5,0)..controls (3.4,1) and (3,1.5)..(4.5,3.5)--(5.5,3.5)..controls (3.2,1.5) and (2.8,1)..(2.5,0);

\filldraw[white,fill=gray!80](2.5,0)..controls (2.8,1) and (3.2,1.5)..(5.5,3.5)--(6.6,3.5)..controls (4,1.5) and (3.8,1)..(3.5,0);
\filldraw[white,fill=gray!80](2.5,0)..controls (2.6,1) and (2.2,1.5)..(3.5,3.5)--(4.5,3.5)..controls (3,1.5) and (3.4,1)..(3.5,0);
\filldraw[white,fill=gray!80](2.5,0)..controls (1.9,2) and (1,2.8)..(0.4,3.5)--(1.3,3.5)..controls (2,2.6) and (3.2,1)..(3.5,0);
\filldraw[white,fill=gray!80](2.5,0)..controls (1,1) and (0.5,1.8)..(-1.8,3.5)--(-0.8,3.5)..controls (1,2) and (2,1)..(3.5,0);
\draw[->](-1.8,0)--(7,0)node[left,below]{$t=0$};
\draw[dashed](-1.8,1.7)--(7,1.7)node[right,below]{$t=t_0^{(\eta)}$};




\draw [color=black](3.5,0)..controls (3.8,1) and (4,1.5)..(6.6,3.5);

\draw [color=black](2.5,0)..controls (2.8,1) and (3.2,1.5)..(5.5,3.5);
\node [below] at(5.3,3.2){$\mathcal{R}_1$};

\draw [color=black](3.5,0)..controls (3.4,1) and (3,1.5)..(4.5,3.5);

\draw [color=black](2.5,0)..controls (2.6,1) and (2.2,1.5)..(3.5,3.5);
\node [below] at(3.55,3.2){$\mathcal{R}_2$};

\draw [color=black](3.5,0)..controls (3.2,1) and (2,2.6)..(1.3,3.5);
\node [below] at(1.4,3.2){$\mathcal{R}_3$};

\draw [color=black] (2.5,0)..controls (1.9,2) and (1,2.8)..(0.4,3.5);

\draw [color=black](3.5,0)..controls (2,1) and (1,2)..(-0.8,3.5);
\node [below] at(-0.5,3.2){$\mathcal{R}_{4}$};

\draw [color=black] (2.5,0)..controls (1,1) and (0.5,1.8)..(-1.8,3.5);

\end{tikzpicture}
\caption{\small{\bf Separation of four characteristic strips.} In this picture, the domain is divided into two regions by the dashed line $t=t_0^{(\eta)}$: (a) The region above the line $t=t_0^{(\eta)}$ denotes the separating domain where the four strips are disjoint. (b) The region under $t=t_0^{(\eta)}$ denotes the non-separating domain where the strips overlap with each other. }
\label{fig1}
\end{figure}
\subsection{$L^\infty$-estimates}
Since the eigenvalues of system \eqref{ei} are distinct, system \eqref{2d2} is a strictly hyperbolic system. Note that all the coefficients of system \eqref{2d2} are uniformly bounded and we assume that they are controlled by a uniform constant $\Gamma$. In this subsection, we prove the $L^\infty$ type estimates for the decomposed system \eqref{eqrho}-\eqref{eqv}. We first derive the a priori estimates. They are deduced in two regions: the non-separating region $[0,t_0^{(\eta)}]$ and the separating region $[t_0^{(\eta)},T]$ as in Figure \ref{fig1}.

{\it Estimates for $[0,t_0^{(\eta)}]$.} Within this time interval, all the characteristic strips might overlap with each other. We first consider
\begin{equation}\label{W}
W(t)=\max_i\sup_{(x',t')\atop 0\leq t'\leq t}\big|w_i(x',t')\big|.
\end{equation}
By \eqref{eqw}, we have
\begin{equation*}
  \frac{\partial}{\partial s_i}|w_i|\leq \Gamma W^2.
\end{equation*}
Via comparing $w_i$ with the solution to the following ODE
\begin{equation*}
  \left\{\begin{array}{ll}
  \frac{d}{dt}Y=\Gamma Y^2,\\
  Y(0)=W_0^{(\eta)},
  \end{array}
  \right.
\end{equation*}
we obtain
\begin{equation}\label{wt0}
  |w_i|\leq Y(t)=\frac{W_0^{(\eta)}}{1-\Gamma W_0^{(\eta)} t}\quad \text{for}\quad t<\min\Big\{\frac1{\Gamma W_0^{(\eta)}},t_0^{(\eta)}\Big\}.
\end{equation}
Since $ t_0^{(\eta)}=O(\eta)$, it holds that
\begin{equation}\label{gwt}
  \Gamma W_0^{(\eta)} t_0^{(\eta)}=O(\eta W_0^{(\eta)})=O(\theta).
\end{equation}
Applying \eqref{gwt} to \eqref{wt0}, for a small parameter $\varepsilon\in(0,\frac1{100}]$, we deduce that
\begin{equation*}
  |w_i(x,t)|\leq (1+\varepsilon)W_0^{(\eta)}\quad \text{for any}\ x\in\mathbb{R}\ \text{and}\ t\in[0,t_0^{(\eta)}].
\end{equation*}
This implies that
\begin{equation}\label{y15}
  |W(t)|\leq (1+\varepsilon)W_0^{(\eta)}\quad \text{for}\  t\in[0,t_0^{(\eta)}].
\end{equation}

Now we consider the exterior region of the characteristic strips. Using characteristic coordinates, any point outside of the characteristic strip $\mathcal{R}_i$ can be labelled by $(z_i',s_i')$ with $z'_i\notin[\eta,2\eta]$, and it obviously holds $w_i^{(\eta)}(z_i',0)=0$. Integrating  \eqref{eqw} along the characteristic $\mathcal{C}_i$, we obtain
\begin{equation}
  V(t)=O(\int_0^{t_0^{(\eta)}}w_iw_jds_i)=O(\eta [W(t)]^2)=O(\eta [W_0^{(\eta)}]^2).
\end{equation}

To bound $S$, we consider the equation of the inverse density $\rho_i$, and we have
\begin{equation}
  \frac{\partial\rho_i}{\partial s_i}=O(\rho_iW).
\end{equation}
Integrating the above equation along the characteristic $\mathcal{C}_i$, we obtain
\begin{equation}\label{rho i}
  \rho_i(z_i,t)=\rho_i(z_i,0)\exp\big(O(tW(t))\big).
\end{equation}
Note that by the definitions of \eqref{flow} and \eqref{dense} we have
\begin{equation}\label{319}
  \rho_i(z_i,0)=1.
\end{equation}
Then it follows from \eqref{y15} and (\ref{rho i}) that
\begin{equation}\label{320}
  \rho_i(z_i,t)=\exp\big(O(\eta W_0^{(\eta)})\big)>0\ \text{for any}\ t\in[0,t_0^{(\eta)}].
\end{equation}
For a fixed $\eta$, we can choose sufficiently small $\theta$ such that
\begin{equation}
 1-\varepsilon\leq\rho_i(z_i,t)\leq1+\varepsilon\ \text{for any}\ t\in[0,t_0^{(\eta)}].
\end{equation}
Hence we have
\begin{equation}\label{y16}
  S(t)=O(1)\ \text{for any}\ t\in[0,t_0^{(\eta)}].
\end{equation}

Next, we estimate the supremum of $v_i$. Note that for $v_i$ we have
\begin{equation*}
\frac{\partial v_i}{\partial s_i}=O(S(t)[W(t)]^2).
\end{equation*}
Integrating along $\mathcal{C}_i$ and using \eqref{y15} and \eqref{y16}, for $t\in[0,t_0^{(\eta)}]$, we get
\begin{equation}
  J(t)=O(W_0^{(\eta)}+t[W(t)]^2)=O(W_0^{(\eta)}+\eta [W_0^{(\eta)}]^2)=O(W_0^{(\eta)}).
\end{equation}

Finally, we obtain the estimate for $\bar{U}$ from
\begin{equation}\label{phi}
  \Phi(x,t)=\int_{X_4(\eta,t)}^x\frac{\partial \Phi(x',t)}{\partial x} dx'=\int_{X_4(\eta,t)}^x\sum_{k}w_kr_k(x',t)dx'.
\end{equation}
This indicates that
\begin{equation}
   \bar{U}(t)=\sup_{(x',t')\atop 0\leq t'\leq t}|\Phi(x',t')|=O\big(W(t)(\eta+(\bar{\lambda}_1-\underline{\lambda}_4)t)\big)=O(\eta W_0^{(\eta)})\ \text{for any}\ t\in[0,t_0^{(\eta)}].
\end{equation}

{\it  Estimates for $[t_0^{(\eta)},T]$.} Due to the strict hyperbolicity, the four aforementioned characteristic strips $\mathcal{R}_i$ are well separated after $t_0^{(\eta)}$.

We firstly estimate $S(t)$, i.e., the supremum of the inverse densities. If $(x,t)\in \mathcal{R}_i$, we have
\begin{equation}\label{335}
  \frac{\partial \rho_i}{\partial s_i}=O(J+ V S).
\end{equation}
Thus, by integrating \eqref{335} along the characteristic $\mathcal{C}_i$, we derive that
\begin{equation}
\rho_i(z_i, t)=\rho_0(z_i,0)+\int_0^tO(J+ V S)ds_i.
\end{equation}
This gives
\begin{equation}
  S(t)=O(1+tJ+ tV S).
\end{equation}

We then bound $J(t)$, the supremum of $\{v_i(x,t)\}_{i=1,2,3,4}$ with $(x,t)\in \mathcal{R}_i$. By \eqref{eqv}, we have
\begin{equation}
  \frac{\partial v_i}{\partial s_i}=O(VJ+ V^2 S).
\end{equation}
Via integration, we get
\begin{equation}\label{548}
  J(t)=O(W_0^{(\eta)}+tVJ+ tV^2 S).
\end{equation}

Next we derive the upper bounds of $w_i$ outside the corresponding characteristic strip $\mathcal{R}_i$. From  \eqref{eqw}, we have that the evolution equations satisfy:
\begin{equation}\label{558}
  \frac{\partial w_i}{\partial{s_i}}=O(V^2)+O\Big(\sum_{k\neq i}w_k\Big)V+O\Big(\sum_{m\neq i,k\neq i\atop m\neq k}w_mw_k\Big).
\end{equation}
Note that $\mathcal{C}_i$ starts from $z_i\notin[\eta,2\eta]$ and ends at $(x,t)\notin\mathcal{R}_i$. When $t'\geq t_0^{(\eta)}$, for any point $\big(X_i(z_i,t'),t'\big)\in\mathcal{C}_i$, it holds either $\big(X_i(z_i,t'),t'\big)\in\big(\mathbb{R}\times[t_0^{(\eta)},t]\big)\setminus\bigcup_{k}\mathcal{R}_k$ or $\big(X_i(z_i,t'),t'\big)\in\mathcal{R}_k$ for some $k\neq i$.
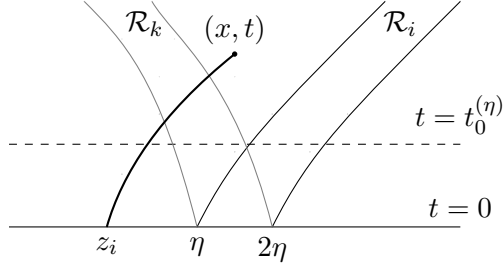
\begin{center}
\begin{tikzpicture}
\draw(0,0)--(6,0)node[left,above]{$t=0$};
\draw[dashed](0,1.1)--(6,1.1)node[left,above]{$t=t_0^{(\eta)}$};
\node [below]at(3.5,0){$2\eta$};
\node [below]at(2.5,0){$\eta$};

\filldraw [black] (3.5,0) circle [radius=0.01pt]
(4,1) circle [radius=0.01pt]
(4.5,1.5) circle [radius=0.01pt]
(6,3) circle [radius=0.01pt];
\draw [color=black](3.5,0)..controls (4,1) and (4.5,1.5)..(6,3);

\filldraw [black] (2.5,0) circle [radius=0.01pt]
(3,1) circle [radius=0.01pt]
(3.5,1.5) circle [radius=0.01pt]
(5,3) circle [radius=0.01pt];
\draw [color=black](2.5,0)..controls (3,1) and (3.5,1.5)..(5,3);
\node [below] at(5.2,3){$\mathcal{R}_{i}$};

\filldraw [gray] (3.5,0) circle [radius=0.01pt]
(3,2) circle [radius=0.01pt]
(2.2,2.5) circle [radius=0.01pt]
(1.9,3) circle [radius=0.01pt];
\draw [color=gray](3.5,0)..controls (3,2) and (2.2,2.5)..(1.9,3);
\node [below] at(1.8,3){$\mathcal{R}_{k}$};

\filldraw [gray] (2.5,0) circle [radius=0.01pt]
(2,2) circle [radius=0.01pt]
(1.5,2.5) circle [radius=0.01pt]
(1,3) circle [radius=0.01pt];
\draw [color=gray] (2.5,0)..controls (2,2) and (1.5,2.5)..(1,3);

\filldraw [black] (1.3,0) circle [radius=0.01pt]
(1.5,0.7) circle [radius=0.01pt]
(2,1.5) circle [radius=0.01pt]
(3,2.3) circle [radius=0.8pt];
\draw [color=black,thick](1.3,0)..controls (1.5,0.7) and (2,1.5)..(3,2.3);
\node[above]at (3,2.3){$(x,t)$};
\node[below]at (1.3,0){$z_i$};

\end{tikzpicture}
\end{center}
For the term $O\Big(\sum_{m\neq i,k\neq i\atop m\neq k}w_mw_k\Big)$ on the RHS of \eqref{558}, if $(x,t)\notin\mathcal{R}_i$, then there are only three scenarios: $(x,t)$ stays in $\mathcal{R}_m$, or $(x,t)$ stays in $\mathcal{R}_k$, or $(x,t)$ stays out of all the characteristics. In all of these three cases, the term $O\Big(\sum_{m\neq i,k\neq i\atop m\neq k}w_mw_k\Big)$ can be absorbed by $O\Big(\sum_{k\neq i}w_k\Big)V$, which is the second term on the RHS of \eqref{558}. Denote $I_k^i=\{t'\in[t_0^{(\eta)},t]:(x,t')\in\mathcal{C}_i\bigcap\mathcal{R}_k\}$  for $k\neq i$. Integrating \eqref{558} along $\mathcal{C}_i$ and using $w_i^{(\eta)}(z_i,0)=0$, for $(x,t)\notin\mathcal{R}_i$ we have that
\begin{equation}\label{559}
\begin{split}
    w_i(x,t)=&O\Big(tV^2+V\sum_{k\neq i}\int_0^tw_k\big(X_i(z_i,t'),t'\big)dt'\Big)\\
    =&O\Big(tV^2+V\sum_{k\neq i}\int_0^{t_0^{(\eta)}}w_k\big(X_i(z_i,t'),t'\big)dt'\Big)\\
    &+O\Big(V\sum_{k\neq i}\int_{t_0^{(\eta)}}^t w_k\big(X_i(z_i,t'),t'\big)dt'\Big)\\
    =&O\Big(tV^2+\eta[W_0^{(\eta)}]^2+V\sum_{k\neq i}\underbrace{\int_{I_k^i}w_k\big(X_i(z_i,t'),t'\big)dt'}_{M}\Big).
\end{split}
\end{equation}
Here  we use the fact $V(t)\leq W(t)=O(W_0^{(\eta)})$ for $t\leq t_0^{(\eta)}$. Then we employ the bi-charateristic coordinates to bound $M$. And we have
\begin{equation}\label{560}
  \begin{split}
    &\int_{I_k^i}w_k\big(X_i(z_i,t'),t'\big)dt'\\
    =&O\Big(\int_{y_k\in[\eta,2\eta]}\Big|\frac{\rho_k\big(y_k,t'(y_i,y_k)\big)}{\lambda_i-\lambda_k}w_k\big(y_k,t'(y_i,y_k)\big)\Big|dy_k\Big)\\
    =&O(\eta J).
  \end{split}
\end{equation}
By \eqref{559} and \eqref{560}, we obtain
\begin{equation}\label{355}
  V(t)=O(tV^2+\eta [W_0^{(\eta)}]^2+\eta VJ).
\end{equation}

The final step is to estimate $|\Phi|$. If $(x,t)$ does not belong to any characteristic strip, we go back to \eqref{phi} and obtain that
\begin{equation}\label{u1}
 \bar{U}(t)=O\big((\eta+(\bar{\lambda}_1-\underline{\lambda}_4)t)V\big).
\end{equation}
If $(x,t)\in\mathcal{R}_k$ for some $k$, using characteristic coordinates, we have
\begin{equation}\label{u2}
\begin{split}
  |\Phi(x,t)|=&\Big|\int_{X_k(\eta,t)}^x\frac{\partial \Phi(x',t)}{\partial x} dx'\Big|=\Big|\int_{X_k(\eta,t)}^x\sum_{k}w_kr_k(x',t)dx'\Big|\\
  =&O\Big(\int_{X_k(\eta,t)}^{X_k(2\eta,t)}|w_k(x',t)|dx'\Big)=O\Big(\int_\eta^{2\eta}|w_k(x',t)|\rho_kdz_k\Big)\\
  =&O(\eta J).
  \end{split}
\end{equation}
Combining \eqref{u1} and \eqref{u2}, we get
\begin{equation}
   \bar{U}(t)=O(\eta J+\eta V+\eta tV).
\end{equation}

In summary, if $\Phi\in C^2(\mathbb{R}\times[0,T],B_{2\kappa}^4(0))$ is a solution to \eqref{2d2} for some $T>0$, we then obtain the following a priori estimates:
\begin{align*}
  S&=O(1+tJ+tVS), &J&=O(W_0^{(\eta)}+t VJ+tV^2 S),\\
 V&=O\Big(\eta [W_0^{(\eta)}]^2+tV^2+\eta VJ\Big),&\bar{U}&=O(\eta J+\eta V+\eta tV),
\end{align*}
where $t<T$.

With a bootstrap argument as in Christodoulou-Perez \cite{christodoulou}, we further derive the following bounds:
\begin{align*}
  S(t)&=O(1), &J(t)&=O(W_0^{(\eta)}),\\\
 V(t)&=O\Big(\eta [W_0^{(\eta)}]^2\Big),&\bar{U}(t)&=O(\eta W_0^{(\eta)}),
\end{align*}
for any $t\in[0,T_{\eta}^*)$ and
\begin{equation}
T_{\eta}^*\leq O\Big(\frac{1}{W_0^{(\eta)}}\Big).
\end{equation}
In particular, we make the following bootstrap assumptions
\begin{equation}
 tV\leq \theta^{\frac12}, \qquad J\leq \theta^{\frac12}.
\end{equation}
Employing the a priori estimates, it holds that
\begin{equation*}
\begin{split}
S(t)=O(1+tJ+\theta^{\frac12}S)&\qquad\Rightarrow\qquad S(t)=O(1+tJ),\\
 V(t)=O\Big(\eta [W_0^{(\eta)}]^2+\theta^{\frac12}V+\eta  \theta^{\frac12}V\Big)&\qquad\Rightarrow\qquad V(t)=O\Big(\eta [W_0^{(\eta)}]^2\Big),
\end{split}
\end{equation*}
and
\begin{equation}
\begin{split}
J(t)=O(W_0^{(\eta)}+\theta^{\frac12}J+\theta^{\frac12}VS)\qquad&\Rightarrow\qquad J(t)=O\big(W_0^{(\eta)}+\theta^{\frac12}V(1+tJ)\big)\\
&\Rightarrow\qquad J(t)=O\big(W_0^{(\eta)}+\theta^{\frac12}V\big)=O\big(W_0^{(\eta)}\big).
\end{split}
\end{equation}
Then for $t< O(\frac{1}{W_0^{(\eta)}})$, we prove $S(t)=O(1)$ and $\bar{U}(t)=O(\eta W_0^{(\eta)})$. Choosing $\theta$ to be sufficiently small, we hence prove $\Phi\in B_\kappa^4(0)$.

\subsection{The first shock forms in $\mathcal{R}_1$ at $t=T_{\eta}^*$}
In the following, we deduce lower bound and upper bound for $T_{\eta}^*$. When $t$ goes to $T_{\eta}^*$, a shock forms in $\mathcal{R}_1$. When the shock forms, the inverse density $\rho_1$ becomes $0$. With the equation for $\rho_1$:
\begin{equation*}
  \frac{\partial\rho_1}{\partial s_1}=c_{11}^1(\Phi)v_1+O\Big(\sum_{ k\neq 1}w_k\Big)\rho_1,
\end{equation*}
and the fact $c_{11}^1(\Phi)<0$, we have
\begin{equation}\label{381}
  -|c_{11}^1||v_1|-\Big|O\Big(\sum_{ k\neq 1}w_k\Big)\Big|\rho_1\leq\frac{\partial\rho_1}{\partial s_1}\leq-|c_{11}^1||v_1|+\Big|O\Big(\sum_{ k\neq 1}w_k\Big)\Big|\rho_1.
\end{equation}
Since $|\Phi|=O(\eta W_0^{(\eta)})\leq\kappa$, we can choose $\theta$ sufficiently small such that
\begin{equation}\label{y37}
 (1-\varepsilon)|c_{11}^1(0)|\leq |c_{11}^1(\Phi)|\leq (1+\varepsilon)|c_{11}^1(0)|.
\end{equation}
Using bi-characteristic coordinates, we derive
\begin{equation*}
  \int_0^{t}\sum_{k\neq 1}w_k(X_1(z_1,t'),t')dt'=O(\eta W_0^{(\eta)}+\eta J)=O(\eta W_0^{(\eta)}).
\end{equation*}
Therefore, we obtain
\begin{equation}\label{385}
  1-\varepsilon\leq \exp{\Big(\int_0^{t}O\big(\sum_{k\neq 1}w_k(X_1(z_1,t'),t')\big)dt'\Big)}\leq1+\varepsilon,
\end{equation}
and
\begin{equation}\label{386}
  1-\varepsilon\leq \exp{\Big(-\int_0^{t}O\big(\sum_{k\neq 1}w_k(X_1(z_1,t'),t')\big)dt'\Big)}\leq1+\varepsilon.
\end{equation}
Applying Gr\"{o}nwall's inequality to \eqref{381}, together with \eqref{y37}-\eqref{386}, we get
\begin{equation}\label{y38}
\begin{split}
  (1-\varepsilon)\Big(1-(1+\varepsilon)^2|c_{11}^1(0)|\int_0^t|v_1(z,t')|dt'\Big) \leq&\rho_1(z,t)\\
\leq& (1+\varepsilon)\Big(1-(1-\varepsilon)^2|c_{11}^1(0)|\int_0^t|v_1(z,t')|dt'\Big).
\end{split}
\end{equation}
Moreover, for $v_1$, we integrate
\begin{equation}\label{y39}
\frac{\partial v_1}{\partial s_1}=O\Big(\sum_{m\neq 1}w_m\Big)v_1+O\Big(\sum_{m,k\neq1}w_mw_k\Big)\rho_1
\end{equation}
along $\mathcal{C}_1$ and obtain
\begin{equation}\label{79}
\begin{split}
  v_1(z,t)\leq w_1^{(\eta)}(z,0)+O(tVJ+tV^2S)
  =w_1^{(\eta)}(z,0)+O(\eta [W_0^{(\eta)}]^2).
\end{split}\end{equation}
Noting $W_0^{(\eta)}=w_1^{(\eta)}(z_0,0)$, the above inequality yields
\begin{equation*}
  v_1(z_0,t)\leq (1+\varepsilon) W_0^{(\eta)}.
\end{equation*}
By using the first inequality of \eqref{y38}, we arrive at
\begin{equation}\label{711}
  \rho_1(z_0,t)\geq(1-\varepsilon)\Big(1-(1+\varepsilon)^3|c_{11}^1(0)|tW_0^{(\eta)}\Big).
\end{equation}
This shows that $\rho_1(z_0,t)>0$ when $t<\frac1{(1+\varepsilon)^3|c_{11}^1(0)|W_0^{(\eta)}}$.

On the other hand, noticing that
\begin{equation*}
  O\Big(\sum_{m\neq 1}\int_0^tw_mdt'\Big)=O(\eta W_0^{(\eta)}+\eta J)=O(\eta W_0^{(\eta)})
\end{equation*}
and
\begin{equation*}
\begin{split}
  O\Big(\sum_{m,k\neq1}\int_0^tw_mw_kdt'\Big)
  =O(t V^2)
  =O(\eta^2 [W_0^{(\eta)}]^3),
  \end{split}
\end{equation*}
via Gr\"{o}nwall's inequality, we have that at $z=z_0$ it holds
\begin{equation} \label{vlow}
 v_1(z_0,t)\geq(1-\varepsilon)[W_0^{(\eta)}-\varepsilon W_0^{(\eta)}]=(1-\varepsilon)^2W_0^{(\eta)}.
\end{equation}
By the second inequality of \eqref{y38}, we get
\begin{equation*}
  \rho_1(z_0,t)\leq (1+\varepsilon)\Big(1-(1-\varepsilon)^4|c_{11}^1(0)|tW_0^{(\eta)}\Big).
\end{equation*}
Together with \eqref{711} and \eqref{vlow}, we conclude that there exists a finite $T_\eta^*$ (shock formation time) such that
\begin{equation*}
  \lim_{t\rightarrow T_\eta^*} \rho_1(z_0,t)=0, \quad \text{and} \quad \lim_{t\rightarrow T_\eta^*} w_1(z_0,t)=+\infty.
\end{equation*}
And $T_\eta^*$ obeys
\begin{equation} \label{Tshock}
 \frac1{(1+\varepsilon)^3|c_{11}^1(0)|W_0^{(\eta)}}\leq T_\eta^*\leq\frac{1}{(1-\varepsilon)^4|c_{11}^1(0)|W_0^{(\eta)}}.
\end{equation}

\subsection{No other shock forms before $T_\eta^*$}
In the above, we prove $w_1$ tends to infinity when $t$ goes to $T_\eta^*$. In this subsection we further prove that $\Big\{\bar W_i=\sup_{(z_i,s_i)\atop z_i\in[\eta,2\eta],\ 0\leq s_i\leq t}w_i(z_i,s_i)\Big\}_{i=2,3,4}$ are all bounded when $t\leq T_\eta^*$.

Invoking the estimate of $V$, for $(x,t)\in\mathcal{R}_i$, it holds that
\begin{equation}\label{eqw7}
\begin{split}
  \frac{\partial w_i}{\partial{s_i}}=&-c^i_{ii}w_i^2+O\Big(\sum_{k\neq i}w_k\Big)w_i+O\Big(\sum_{m\neq i,k\neq i\atop m\neq k}w_mw_k\Big)\\
  =&-c^i_{ii}w_i^2+O(V)w_i+O(V^2)\\
  =&-c^i_{ii}w_i^2+O\big(\eta[W_0^{(\eta)}]^2\big)w_i+O\big(\eta^2[W_0^{(\eta)}]^4\big).
  \end{split}
\end{equation}
Equivalently, we have
\begin{small}
\begin{equation} \label{w7ode}
\frac{d}{ds}\Big[\exp\Big(O\big(\eta[W_0^{(\eta)}]^2\big)s\Big)w_i\Big]=-\exp\Big(O\big(\eta[W_0^{(\eta)}]^2\big)s\Big)c^i_{ii}w_i^2+\exp\Big(O\big(\eta[W_0^{(\eta)}]^2\big)s\Big)O\big(\eta^2[W_0^{(\eta)}]^4\big).
\end{equation}
\end{small}
Since
\begin{equation}
\exp\Big(O\big(\eta[W_0^{(\eta)}]^2\big)s\Big)\leq \exp\Big(O\big(\eta[W_0^{(\eta)}]^2\big)T_\eta^*\Big)\leq \exp\Big(O\big(\eta W_0^{(\eta)}\big)\Big)=O(e^\theta),
\end{equation}
by choosing $\theta$ sufficiently small, we have
\begin{equation} \label{exp}
1-\varepsilon\leq\exp\Big(O\big(\eta[W_0^{(\eta)}]^2\big)s\Big)\leq 1+\varepsilon.
\end{equation}
Thus, invoking \eqref{exp} in \eqref{w7ode}, we obtain
\begin{equation}\label{deqw7}
\begin{split}
\frac{d}{ds}\Big[\exp\Big(O\big(\eta[W_0^{(\eta)}]^2\big)s\Big)w_i\Big]
\leq Cw_i^2+O\big(\eta^2[W_0^{(\eta)}]^4\big).
\end{split}
\end{equation}
Integrating \eqref{deqw7} along $\mathcal{C}_i$, we deduce that
\begin{equation}\label{621}
\begin{split}
\bar W_i&\leq O\big(w_i(z,0)+t\bar W_i^2+t\eta^2[W_0^{(\eta)}]^4\big)\\
&\leq O\big([W_0^{(\eta)}]^2+t\bar W_i^2+\eta^2[W_0^{(\eta)}]^3\big)\\
&\leq O\big([W_0^{(\eta)}]^2+t\bar W_i^2\big).
\end{split}
\end{equation}
Then we introduce an additional bootstrap assumption
\begin{equation} \label{boottcw}
t\bar W_i\leq \theta^{\frac12}.
\end{equation}
By \eqref{621}, it holds that
\begin{equation} \label{Wcc7}
\bar W_i\leq  O\big([W_0^{(\eta)}]^2\big)
\end{equation}
and the bootstrap assumption \eqref{boottcw} can be improved to
\begin{equation}
t\bar W_i\leq  O(W_0^{(\eta)})=O(\theta)< \theta^{\frac12}.
\end{equation}
This implies the boundedness of $\bar W_i$.

\section{$H^\frac{11}{4}$ ill-posedness} \label{pfill}
Now we prove the $H^\frac{11}{4}$ ill-posedness stated in Theorem \ref{3D}. First note that via our construction \eqref{dataw0} and \eqref{data2} for initial data, using the estimates \eqref{Tshock} for the shock formation time $T^*_\eta$, it follows immediately that
\begin{equation*}
T^*_\eta \to 0 \quad \text{as}\quad \eta \to 0.
\end{equation*}
In the below, we complete the proof of Theorem \ref{3D} by further showing that the solution's $H^2$ norm blows up at the shock formation time $T_{\eta}^*$.

For the initial data constructed in Section \ref{2d data}, without loss of generality, we restrict our attention to the following region
\begin{equation*}
\Omega_0=\{(x,y_2):\frac{{c_2^2(x-\frac{3\eta}{2})}^2}{c_1^2}+(y_2)^2 \leq \frac{c_2^2}{c_1^2} {\left(\frac{\eta}{2} \right)}^2\}.
\end{equation*}
One can verify that the initial data is planar symmetric within $\Omega_0$, i.e., $\hat{w}_i^{(\eta)}(x,x_2)|_{\Omega_0}=\hat{w}_i^{(\eta)}(x)$ for $i=1,2,3,4$. Now we write $\Omega_{T_\eta^*}$ to denote the $T_\eta^*$-slice of the domain of future dependence of $\Omega_0$. Then, for $\Phi\in B_\kappa^4(0)$, $\Omega_{T_\eta^*}$ is uniformly close to an ellipse in 2D defined by
\begin{equation*}
\Omega_{T_\eta^*}\approx\{(x,y_2):\frac{{c_2^2(x-c_1T_\eta^*-\frac{3\eta}{2})}^2}{c_1^2}+(y_2)^2 \leq \frac{c_2^2}{c_1^2} {\left(\frac{\eta}{2} \right)}^2\}.
\end{equation*}
And we have
\begin{equation}
\int_{(x,y_2)\in\Omega_{T_\eta^*}} dy_2 \approx(1+O(\varepsilon))\sqrt{\big(z-\eta+O(\varepsilon)\eta\big)\big(2\eta-z+O(\varepsilon)\eta\big)}.
\end{equation}

To achieve our goal, we first prove that $|\partial_{z_1}\rho_1(z_1,s_1)|$ is bounded.
\begin{prop} \label{dzrho1}
For any $s_1<T_\eta^*$, we have $|\partial_{z_1}\rho_1(z_1,s_1)|\leq C$. Here $C$ is a uniform constant depending only on $\varepsilon,\theta$ and $\eta$.
\end{prop}

\begin{proof}
The proof of this proposition is similar to our proof for the 3D case in \cite{an,an2}. We present the main ideas here. Using \eqref{biytoz}, we get
\begin{equation}\label{rhod}
 \partial_{z_1}\rho_1=\partial_{y_1}\rho_1+\frac{\rho_1}{2\lambda_1}\partial_{s_1}\rho_1.
\end{equation}
To bound $\partial_{z_1}\rho_1$, we start with controlling $\partial_{y_1}\rho_1$. Let
\begin{equation*}
\tau_{1} ^{(4)}:=\partial_{y_1}\rho_1,\quad\pi_{1}^{(4)}:=\partial_{y_1}v_1.
\end{equation*}
Since
\begin{equation}\label{96}
  \partial_{y_4}=\frac{\rho_4}{\lambda_1-\lambda_4}\partial_{s_1}=\frac{\rho_4}{2\lambda_1}\partial_{s_1},
\end{equation}
we have that $\tau_{1}^{(4)}$ obeys
\begin{equation}\label{97}
\begin{split}
  \partial_{y_4}\tau_{1}^{(4)}=&\partial_{y_1}\partial_{y_4}\rho_1=\partial_{y_1}\Big(\frac{\rho_4}{2\lambda_1}\partial_{s_1}\rho_1\Big)=\partial_{y_1}\Big(\frac{\rho_1\rho_4}{2\lambda_1}\sum_m c_{1m}^1w_m\Big)\\
  =&\frac{\rho_4}{2\lambda_1}\Big(c_{11}^1\partial_{y_1}v_1+\sum_{m\neq1}c_{1m}^1w_m \partial_{y_1}\rho_1\Big)+\frac{\rho_1}{2\lambda_1}\Big(\sum_{m\neq4}c_{1m}^1w_m\partial_{y_1}\rho_4+c_{14}^1\partial_{y_1}v_4\Big)\\
  &+\frac{\rho_1\rho_4}{2\lambda_1}\Big(\sum_{m=2,3}c_{1m}^1\frac1{\rho_m}\partial_{y_1}v_m-\sum_{m=2,3}c_{1m}^1\frac{w_m}{\rho_m}\partial_{y_1}\rho_m\Big)\\
  &-\frac{\partial_{y_1}\lambda_1}{2\lambda_1^2}\rho_1\rho_4\sum_m c_{1m}^1w_m+\frac{\rho_1\rho_4}{2\lambda_1}\sum_m \partial_{y_1}c_{1m}^1w_m.
\end{split}
\end{equation}
We still need to estimate $\partial_{y_1}\lambda_1,\partial_{y_1}c_{1m}^1,\partial_{y_1}\rho_m$ and $\partial_{y_1}v_m$ in \eqref{97}. For $\lambda_1$, under bi-characteristic coordinates $(y_1,y_i)$ ($i\neq1$), its derivative satisfies
\begin{equation}\label{9.11}
  \begin{split}
    \partial_{y_1}\lambda_1=&\nabla_\Phi\lambda_1\cdot\partial_{y_1}\Phi=\nabla_\Phi\lambda_1\cdot[\partial_{s_i}X_i\partial_{y_1}t'\partial_{x}\Phi+\partial_{y_1}t'\partial_{t}\Phi]\\
    =&\nabla_\Phi\lambda_1\cdot\Big[\lambda_i\frac{\rho_1}{\lambda_i-\lambda_1}\sum_kw_kr_k+\frac{\rho_1}{\lambda_i-\lambda_1}(-A(\Phi)\sum_kw_kr_k)\Big]\\
    =&O(v_1+\rho_1\sum_{k\neq1}w_k).
  \end{split}
\end{equation}
Similarly, for $c^1_{1m}$, it also holds
\begin{equation}\label{9.12}
  \begin{split}
    \partial_{y_1}c_{1m}^1=&\nabla_\Phi c_{1m}^1\cdot\partial_{y_1}\Phi=\nabla_\Phi c_{1m}^1\cdot[\partial_{s_i}X_i\partial_{y_1}t'\partial_{x}\Phi+\partial_{y_1}t'\partial_{t}\Phi]\\
    =&\nabla_\Phi c_{1m}^1\cdot\Big[\lambda_i\frac{\rho_1}{\lambda_i-\lambda_1}\sum_kw_kr_k+\frac{\rho_1}{\lambda_i-\lambda_1}(-A(\Phi)\sum_kw_kr_k)\Big]\\
    =&O\Big(v_1+\rho_1\sum_{k\neq1}w_k\Big).
  \end{split}
\end{equation}
By \eqref{biytoz} and \eqref{eqrho}-\eqref{eqv}, we get
\begin{equation}\label{9.14}
  \partial_{y_1}\rho_m=\frac{\rho_1}{\lambda_m-\lambda_1}\partial_{s_m}\rho_m=O\Big(\rho_mv_1+\rho_1\rho_m\sum_{k\neq 1}w_k\Big) \quad \text{when}\quad m\neq 1,
\end{equation}
and
\begin{equation}\label{9.15}
  \partial_{y_1}v_m=\frac{\rho_1}{\lambda_m-\lambda_1}\partial_{s_m}v_m=O\Big(\rho_mv_1\sum_{k\neq1}w_k+\rho_1\rho_m\sum_{j\neq1,k\neq1\atop j\neq k}w_jw_k\Big)\quad \text{when}\quad m\neq 1.
\end{equation}
Note that all the coefficients are bounded, it follows from \eqref{9.11}-\eqref{9.15} that \eqref{97} could be rewritten as
\begin{equation}\label{linear}
\begin{split}
  \partial_{y_4}\tau_{1}^{(4)}:=B_{11}\tau_1^{(4)}+B_{12}\pi_1^{(4)}+B_{13},
\end{split}
\end{equation}
where $B_{11},B_{12},B_{13}$ are uniformly bounded constants depending on $\eta$.

In the same fashion, we also obtain the evolution equation for $\pi_1^{(4)}$:
\begin{equation}\label{917}
  \begin{split}
  \partial_{y_4}\pi_{1}^{(4)}
  =&\frac{\rho_4 }{2\lambda_1}\Big(\sum_{p\neq 1,q\neq 1\atop p\neq q}\gamma_{pq}^1w_pw_q \tau_{1}^{(4)}+\sum_{p\neq 1}\gamma_{1p}^1w_p\pi_{1}^{(4)}\Big)\\
  &-\frac{\partial_{y_1}\lambda_1}{4\lambda_1^2}\Big(\sum_{p\neq 1}\gamma_{1p}^1 w_p v_1\rho_4+\sum_{p\neq 1,q\neq 1\atop p\neq q}\gamma_{pq}^1w_pw_q\rho_4\rho_1\Big)\\
    &+\frac{\rho_4\rho_1}{2\lambda_1}\Big(\sum_{p\neq 1}\partial_{y_1}\gamma_{1p}^1 w_p w_m+\sum_{p\neq 1,q\neq 1\atop p\neq q}\partial_{y_1}\gamma_{pq}^1w_pw_q\Big)\\
    &+\frac{\rho_1}{2\lambda_1}\Big(\sum_{p\neq 1}\gamma_{1p}^1 w_p w_1+\sum_{p\neq 1,q\neq 1\atop p\neq q}\gamma_{pq}^1w_pw_q\Big)\partial_{y_1}\rho_4\\
    &+\frac{\rho_4\rho_1}{2\lambda_1}\Big(\sum_{p=2,3}\gamma_{1p}^1\frac{w_1}{\rho_p}+\sum_{p\neq 1,q\neq 1\atop p\neq q}\gamma_{pq}^1\frac{w_q}{\rho_p}\Big)(\partial_{y_1}v_p-w_p\partial_{y_1}\rho_p)\\
    &+\frac{\rho_4\rho_1}{2\lambda_1 }\sum_{p\neq 1,q\neq 1\atop p\neq q}\gamma_{pq}^1\frac{w_p}{\rho_q}(\partial_{y_1}v_q-w_q\partial_{y_1}\rho_q)
    +\frac{\rho_1}{2\lambda_1 }\gamma_{14}^1 w_1\partial_{y_1}v_4\\
    :=&B_{21}\tau_1^{(4)}+B_{22}\pi_1^{(4)}+B_{23}.
  \end{split}
\end{equation}
As in \eqref{linear}, here one can also check that $B_{21},B_{22},B_{23}$ are all uniformly bounded.

Next, we claim the initial data for $\tau_1^{(4)}$ and $\pi_1^{(4)}$ are also finite. Since $\rho_1(z_1,0)=1$, by \eqref{biytoz}, we have that
\begin{equation*}
  \begin{split}
    \tau_1^{(4)}(z_1,0)=&\partial_{z_1}\rho_1(z_1,0)-\frac{\rho_1(z_1,0)}{2\lambda_1}\partial_{s_1}\rho_1(z_1,0)\\
     =&-\frac{1}{2\lambda_1}\sum_{k}c_{1k}^1 w_k(z_1,0)=O(W_0^{(\eta)})<+\infty.
  \end{split}
\end{equation*}
And noticing that $v_1(z_1,0)=w_1^{(\eta)}(z_1,0)$, we similarly deduce that
\begin{equation*}
  \begin{split}
    \pi_1^{(4)}(z_1,0)=&\partial_{z_1}v_1(z_1,0)-\frac{\rho_1(z_1,0)}{2\lambda_1}\partial_{s_1}v_1(z_1,0)\\
    =&\partial_{z_1}w_1(z_1,0)-\frac{1}{2\lambda_1}\sum_{q\neq 1,q\neq p}\gamma_{pq}^1w_p(z_1,0)w_q(z_1,0)\rho_1(z_1,0)\\
    =&O(\partial_{z_1}w_1(z_1,0)+[W_0^{(\eta)}]^2)<+\infty.
  \end{split}
\end{equation*}
Applying Gr\"onwall's inequality to \eqref{linear} and \eqref{917}, for $s_1\leq T_\eta^*$ we obtain that $\tau_{1}^{(4)}:=\partial_{y_1}\rho_1(z_1,s_1)$ is bounded on $z_1\in[\eta,2\eta]$.

Combining all these estimates in \eqref{rhod}, we hence prove
\begin{equation*}
  \partial_{z_1}\rho_1=\partial_{y_1}\rho_1+O(v_1+\sum_{m\neq 1}w_m\rho_1).
\end{equation*}
With the bounds for $J(t)$, $S(t)$ and $V(t)$, which are obtained in Section \ref{aprior}, we conclude that $\partial_{z_1}\rho_1$ is uniformly bounded. \end{proof}

Now we proceed to estimate the $H^2$ norm of the solutions. We have that
\begin{equation*}
\begin{split}
  &\|w_1(\cdot,T_\eta^*)\|_{L^2(\Omega_{T_\eta^*})}^2\\
  \geq &C\int_{\eta}^{2\eta}\Big|\frac{v_1}{\rho_1}(z,T_\eta^*)\Big|^2\rho_1(z,T_\eta^*) \sqrt{\big(z-\eta+O(\varepsilon)\eta\big)\big(2\eta-z+O(\varepsilon)\eta\big)}dz\\
\geq &C{(1-\varepsilon)}^2{(1-2\varepsilon)}^2[W_0^{(\eta)}]^2\int_{z_0}^{z_0^*}\frac{1}{\rho_1(z,T_\eta^*)} \sqrt{\big(z-\eta+O(\varepsilon)\eta\big)\big(2\eta-z+O(\varepsilon)\eta\big)}dz\\
\geq& C\sqrt{(z_0-\eta)(2\eta-z_0^*)}[W_0^{(\eta)}]^2\int_{z_0}^{z_0^*}\frac{1}{\rho_1(z,T_\eta^*)-\rho_1(z_0,T_\eta^*)}dz\\
\geq& C\sqrt{(z_0-\eta)(2\eta-z_0^*)}[W_0^{(\eta)}]^2\int_{z_0}^{z_0^*}\frac{1}{(\sup_{z\in(z_0,z_0^*]}|\partial_z\rho_1|)(z-z_0)}dz\\
=&+\infty.
  \end{split}
\end{equation*}
Here, we crucially use the fact that a shock form at $(z_0,T_\eta^*)$, i.e., $\rho_1(z_0,T_\eta^*)=0.$ Since the region $\Omega_{T_\eta^*}$ stays in $\mathcal{R}_1$, $w_i\big|_{\Omega_{T_\eta^*}}$ are controlled by $V(t)$ as $t\to T_\eta^*$ for $i=2,3,4$. By the boundedness of $\partial_{z_1}\rho_1$ and $|r_k|$, we derive
\begin{equation}
  \begin{split}
    \|\partial_x^2U_1\|_{L^2(\Omega_{T_\eta^*})}=&\|\sum_{k=1}^4w_kr_{k1}\|_{L^2(\Omega_{T_\eta^*})}\geq\|w_1r_{11}\|_{L^2(\Omega_{T_\eta^*})}-\sum_{k=2}^4\|w_kr_{k1}\|_{L^2(\Omega_{T_\eta^*})}\\
    \geq&C \Big[\|w_1\|_{L^2(\Omega_{T_\eta^*})}-\sum_{k=2}^4\|w_k\|_{L^2(\Omega_{T_\eta^*})}\Big]\\
    \geq&C \Big[\|w_1\|_{L^2(\Omega_{T_\eta^*})}-3V(T_\eta^*)|\Omega_{T_\eta^*}|^\frac{1}{2}\Big]\\
    \geq&C\Big[\|w_1\|_{L^2(\Omega_{T_\eta^*})}-3\eta^2[W_0^{(\eta)}]^2\Big].
  \end{split}
\end{equation}
Finally, since $\|w_1\|_{L^2(\Omega_{T_\eta^*})}=+\infty$, we obtain
\begin{equation*}
  \|\partial_x^2U_1\|_{L^2(\Omega_{T_\eta^*})}=+\infty.
\end{equation*}
This concludes the proof of Theorem \ref{3D}.

\end{document}